\documentclass[11pt]{article}
\usepackage{amsfonts}
\usepackage{amsfonts}
\usepackage{amsfonts}
\usepackage{mathrsfs}
\usepackage{bm}
\usepackage{amssymb,amsmath} % font used for R in Real numbers
 \allowdisplaybreaks

\catcode`!=11
\let\!int\int \def\int{\displaystyle\!int}
\let\!lim\lim \def\lim{\displaystyle\!lim}
\let\!sum\sum \def\sum{\displaystyle\!sum}
\let\!sup\sup \def\sup{\displaystyle\!sup}
\let\!inf\inf \def\inf{\displaystyle\!inf}
\let\!cap\cap \def\cap{\displaystyle\!cap}
\let\!max\max \def\max{\displaystyle\!max}
\let\!min\min \def\min{\displaystyle\!min}
\let\!frac\frac \def\frac{\displaystyle\!frac}
\catcode`!=12

\let\oldsection\section
\renewcommand\section{\setcounter{equation}{0}\oldsection}

\allowdisplaybreaks
\def\pf{\it{Proof.}\rm\quad}

\def\Z{\mathbb{Z}}

\newtheorem{thm}{Theorem}[section]

\newtheorem{cor}[thm]{Corollary}

\begin{document}
%%%%%%%%%%%%%%%%%%%% title %%%%%%%%%%%%%%%%%%%%%%%%%%%%%%%%%%%%%%%%%%%%%%%%
\title {\bf Explicit evaluation of harmonic sums}
\author{
{Ce Xu\thanks{Corresponding author. Email: xuce1242063253@163.com} }\\[1mm]
\small School of Mathematical Sciences, Xiamen University\\
\small Xiamen
361005, P.R. China}

\date{}
\maketitle \noindent{\bf Abstract }  In this paper, we obtain some formulae for harmonic sums, alternating harmonic sums and Stirling number sums by using the method of integral representations of series. As applications of these formulae, we give explicit formula of several quadratic and cubic Euler sums through zeta values and linear sums. Furthermore, some relationships between harmonic numbers and Stirling numbers of the first kind are established.
\\[2mm]
\noindent{\bf Keywords} Harmonic number; Euler sum; Riemann zeta function; Stirling number.
\\[2mm]
\noindent{\bf AMS Subject Classifications (2010):} 11M06; 11M32; 33B15

\section{Introduction}
In this paper, the generalized harmonic numbers and alternating harmonic numbers of order $k$ are defined respectively as
 $$\zeta_n(k):=\sum\limits_{j=1}^n\frac {1}{j^k} ,\ L_{n}(k):=\sum\limits_{j=1}^n\frac{(-1)^{j-1}}{j^k},\ 1\leq n, k \in \Z,\eqno(1.1)$$
where $H_n\equiv\zeta_n(1)=\sum\limits_{j=1}^n\frac {1}{j}$ denotes the classical harmonic number.\\
From [12,27], we know that the classical Euler sums are the infinite sums whose general term is product of 1 or $(-1)^{n-1}$, harmonic numbers and alternating harmonic numbers of index $n$ and a power of $n^{-1}$. Therefore, the Euler sums of index
\[{\pi _1}: = \left( {\underbrace {{k_1}, \cdots ,{k_1}}_{{q_1}}, \cdots ,\underbrace {{k_{{m_1}}}, \cdots ,{k_{{m_1}}}}_{{q_{{m_1}}}}} \right),{\pi _2}: = \left( {\underbrace {{h_1}, \cdots ,{h_1}}_{{l_1}}, \cdots ,\underbrace {{h_{{m_2}}}, \cdots ,{h_{{m_2}}}}_{{l_{{m_2}}}}} \right),p\]
are defined by
$${S_{{\pi _1},{\pi _2},p}}\equiv{S_{\prod\limits_{i = 1}^{{m_1}} {k_i^{{q_i}}} ,\prod\limits_{j = 1}^{{m_2}} {h_j^{{l_j}}} ,p}}:=\sum\limits_{n = 1}^\infty  {\frac{{\prod\limits_{i = 1}^{{m_1}} {\zeta _n^{{q_{_i}}}\left( {{k_i}} \right)\prod\limits_{j = 1}^{{m_2}} {L_n^{{l_j}}\left( {{h_j}} \right)} } }}{{{n^p}}}},$$
$${{\bar S}_{{\pi _1},{\pi _2},p}}\equiv{\bar S_{\prod\limits_{i = 1}^{{m_1}} {k_i^{{q_i}}} ,\prod\limits_{j = 1}^{{m_2}} {h_j^{{l_j}}} ,p}}:=\sum\limits_{n = 1}^\infty  {\frac{{\prod\limits_{i = 1}^{{m_1}} {\zeta _n^{{q_{_i}}}\left( {{k_i}} \right)\prod\limits_{j = 1}^{{m_2}} {L_n^{{l_j}}\left( {{h_j}} \right)} {{}{}}} }}{{{n^p}}}}{\left( { - 1} \right)}^{n - 1},$$
where the quantity $\omega  = {\pi _1} + {\pi _2} + p = \sum\limits_{i = 1}^{{m_1}} {\left( {{k_i}{q_i}} \right)}  + \sum\limits_{j = 1}^{{m_2}} {\left( {{h_j}{l_j}} \right)}  + p$ being called the weight and the quantities $m_1,m_2$ being the degree. Here $p (p>1),m_{1},m_{2},q_{i},k_{i},h_{j},l_{j}$ are non-negative integers, $1 \le {k_1} < {k_2} <  \cdots  < {k_{{m_1}}} \in \Z,\;1 \le {h_1} < {h_2} <  \cdots  < {h_{{m_2}}} \in \Z$. For example,
\[{S_{{1^3}{2^5}{3^2},{1^2}{2^3},4}} = \sum\limits_{n = 1}^\infty  {\frac{{H_n^3\zeta _n^5\left( 2 \right)\zeta _n^2\left( 3 \right)L_n^2\left( 1 \right)L_n^3\left( 2 \right)}}{{{n^4}}}},{{\bar S}_{0,{p_1}{p_2},p}} = \sum\limits_{n = 1}^\infty  {\frac{{{L_n}\left( {{p_1}} \right){L_n}\left( {{p_2}} \right)}}{{{n^p}}}{{\left( { - 1} \right)}^{n - 1}}}\]\[{{\bar S}_{{1^2}{2^4}3,{1^2}{2^3},5}} = \sum\limits_{n = 1}^\infty  {\frac{{H_n^2\zeta _n^4\left( 2 \right){\zeta _n}\left( 3 \right)L_n^2\left( 1 \right)L_n^3\left( 2 \right)}}{{{n^5}}}{{\left( { - 1} \right)}^{n - 1}}},{S_{{p_1}{p_2},0,p}} = \sum\limits_{n = 1}^\infty  {\frac{{{\zeta _n}\left( {{p_1}} \right){\zeta _n}\left( {{p_2}} \right)}}{{{n^p}}}} .\]
From the definition of Euler sums, there are altogether four types of linear sums:
\[{S_{p,0,q}} = \sum\limits_{n = 1}^\infty  {\frac{{{\zeta _n}\left( p \right)}}{{{n^q}}}} ,{{\bar S}_{p,0,q}} = \sum\limits_{n = 1}^\infty  {\frac{{{\zeta _n}\left( p \right)}}{{{n^q}}}{{\left( { - 1} \right)}^{n - 1}}} ,{S_{0,p,q}} = \sum\limits_{n = 1}^\infty  {\frac{{{L_n}\left( p \right)}}{{{n^q}}}} ,{{\bar S}_{0,p,q}} = \sum\limits_{n = 1}^\infty  {\frac{{{L_n}\left( p \right)}}{{{n^q}}}{{\left( { - 1} \right)}^{n - 1}}} .\]
The evaluation of linear sums in
terms of values of Riemann zeta function and polylogarithm function at positive integers is known when $(p,q) = (1,3), (2,2)$, or $p + q$ is odd [6,13]. For instance, we have
$${{\bar S}_{1,0,3}}=\sum\limits_{n = 1}^\infty  {\frac{{{H_n}{}}}{{{n^3}}}}{\left( { - 1} \right)}^{n - 1}  =  - 2{\rm Li}{_4}\left( {\frac{1}{2}} \right) + \frac{{11{}}}{{4}} \zeta(4) + \frac{{{1}}}{{2}}\zeta(2){\ln ^2}2 - \frac{1}{{12}}{\ln ^4}2 - \frac{7}{4}\zeta \left( 3 \right)\ln 2,$$
$${{\bar S}_{0,1,3}}=\sum\limits_{n = 1}^\infty  {\frac{{{L_n}\left( 1 \right)}}{{{n^3}}}{{\left( { - 1} \right)}^{n - 1}} = \frac{3}{2}\zeta \left( 4 \right)}  + \frac{1}{2}\zeta \left( 2 \right){\ln ^2}2 - \frac{1}{{12}}{\ln ^4}2 - 2{\rm Li}{_4}\left( {\frac{1}{2}} \right).$$
In [12], Philippe Flajolet and Bruno Salvy gave explicit reductions to zeta values for all linear sums
${S_{p,0,q}},{{\bar S}_{p,0,q}},{S_{0,p,q}},{{\bar S}_{0,p,q}}$ when $p+q$ is an odd weight.
The relationship between the
values of the Riemann zeta function and Euler sums has been studied by many authors, for example see [2-3,6-14,16-27]. The Riemann zeta function and  alternating Riemann zeta function are defined respectively by [1,4,5]
$$\zeta(s):=\sum\limits_{n = 1}^\infty {\frac {1}{n^{s}}},\Re(s)>1, $$
and
$$\bar \zeta \left( s \right) := \sum\limits_{n = 1}^\infty  {\frac{{{{\left( { - 1} \right)}^{n - 1}}}}{{{n^s}}}} ,\;{\mathop{\Re}\nolimits} \left( s \right) \ge 1.$$
The general multiple zeta functions is defined as
\[\zeta \left( {{s_1},{s_2}, \cdots ,{s_m}} \right) := \sum\limits_{{n_1} > {n_2} >  \cdots  > {n_m} > 0}^{} {\frac{1}{{n_1^{{s_1}}n_2^{{s_2}} \cdots n_m^{{s_m}}}}}, \]
where ${s_1} +  \cdots  + {s_m}$ is called the weight and $m$ is the multiplicity.
In this paper, we show that the Euler-type sums with harmonic numbers:
\[\sum\limits_{n = 1}^\infty  {\frac{{{\zeta _n}\left( m \right)}}{{n\left( {n + k} \right)}}} ,\sum\limits_{n = 1}^\infty  {\frac{{{L_n}\left( m \right)}}{{n\left( {n + k} \right)}}},1 \le m,k \in \Z \]
can be expressed in terms of series of Riemann zeta values and harmonic numbers. We also provide an explicit evaluation of \[\left( {p - 1} \right)!\sum\limits_{n = p - 1}^\infty  {\frac{{S\left( {n + 1,p} \right)}}{{n!n\left( {n + k} \right)}}},2 \le p \in \Z,1 \le k \in \Z \]
in a closed form in terms of zeta values and the complete exponential Bell polynomial ${Y_k}\left( n \right)$, ${S\left( {n,k} \right)}$ stands for the Stirling number of the first kind.
Specifically, we investigate closed-form representations for sums of the following form:
\[{S_{{1^2},0,p}} = \sum\limits_{n = 1}^\infty  {\frac{{H_n^2}}{{{n^p}}}} ,{S_{1,1,p}} = \sum\limits_{n = 1}^\infty  {\frac{{{H_n}{L_n}\left( 1 \right)}}{{{n^p}}}} ,{{\bar S}_{1,1,p}} = \sum\limits_{n = 1}^\infty  {\frac{{{H_n}{L_n}\left( 1 \right)}}{{{n^p}}}{{\left( { - 1} \right)}^{n - 1}}} ,\]
\[{{\bar S}_{{1^2},0,p}} = \sum\limits_{n = 1}^\infty  {\frac{{H_n^2}}{{{n^p}}}{{\left( { - 1} \right)}^{n - 1}}} ,{S_{0,{1^2},p}} = \sum\limits_{n = 1}^\infty  {\frac{{L_n^2\left( 1 \right)}}{{{n^p}}}} ,{{\bar S}_{0,{1^2},p}} = \sum\limits_{n = 1}^\infty  {\frac{{L_n^2\left( 1 \right)}}{{{n^p}}}{{\left( { - 1} \right)}^{n - 1}}} .\]
Furthermore, we evaluate several other series involving harmonic numbers. For example
\[\sum\limits_{n = 1}^\infty  {\frac{{H_n^2}}{{n\left( {n + k} \right)}}}  = \frac{1}{k}\left\{ \begin{array}{l}
 3\zeta \left( 3 \right) + \frac{{H_k^3 + 3{H_k}{\zeta _k}\left( 2 \right) + 2{\zeta _k}\left( 3 \right)}}{3} \\
  - \frac{{H_k^2 + {\zeta _k}\left( 2 \right)}}{k} - \sum\limits_{i = 1}^{k - 1} {\frac{{{H_i}}}{{{i^2}}}}  + \zeta \left( 2 \right){H_{k - 1}} \\
 \end{array} \right\}.\tag{1.2}\]
\section{Main Theorems and Proof}
In this section, we use certain integral of polylogarithm function representations to evaluate several series with alternating (or non-alternating) harmonic numbers. The polylogarithm function defined as follows
\[{\rm Li}{_p}\left( x \right) := \sum\limits_{n = 1}^\infty  {\frac{{{x^n}}}{{{n^p}}}}, \Re(p)>1,\ \left| x \right| \le 1 .\]
when $x$ takes $1$ and $ -1$, then the function ${{\rm Li}{_p}\left( {x} \right)}$ are reducible to Riemann zeta function and alternating Riemann zeta function, respectively.
\begin{thm} For positive integers $m$ and $k$, then
\[\sum\limits_{n = 1}^\infty  {\frac{{{\zeta _n}\left( m \right)}}{{n\left( {n + k} \right)}}}  = \frac{1}{k}\left\{ {\zeta \left( {m + 1} \right) + \sum\limits_{j = 1}^{m - 1} {{{\left( { - 1} \right)}^{j - 1}}\zeta \left( {m + 1 - j} \right){\zeta _{k - 1}}\left( j \right)}  + {{\left( { - 1} \right)}^{m - 1}}\sum\limits_{i = 1}^{k - 1} {\frac{{{H_i}}}{{{i^m}}}} } \right\},\tag{2.1}\]
\[\sum\limits_{n = 1}^\infty  {\frac{{{L_n}\left( m \right)}}{{n\left( {n + k} \right)}}}  = \frac{1}{k}\left\{ \begin{array}{l}
 \bar \zeta \left( {m + 1} \right) + \sum\limits_{j = 1}^{m - 1} {{{\left( { - 1} \right)}^{j - 1}}\bar \zeta \left( {m + 1 - j} \right){\zeta _{k - 1}}\left( j \right)}  \\
  + {\left( { - 1} \right)^{m - 1}}\ln 2\left( {{\zeta _{k - 1}}\left( m \right) + {L_{k - 1}}\left( m \right)} \right) + {\left( { - 1} \right)^m}\sum\limits_{i = 1}^{k - 1} {\frac{{{{\left( { - 1} \right)}^{i - 1}}}}{{{i^m}}}{L_i}\left( 1 \right)}  \\
 \end{array} \right\}.\tag{2.2}\]
\end{thm}
\pf By the definition of polylogarithm function and Cauchy product formula, we can verify that
\[\frac{{{\rm Li}{_m}\left( x \right)}}{{1 - x}} = \sum\limits_{n = 1}^\infty  {{\zeta _n}\left( m \right){x^n}}, \ x\in (-1,1),\tag{2.3}\]
\[- \frac{{{\rm Li}{_m}\left( { - x} \right)}}{{1 - x}}=\sum\limits_{n = 1}^\infty  {{L_n}\left( m \right){x^n}},\ x\in (-1,1).\tag{2.4}\]
Multiplying (2.3) and (2.4) by $x^{r-1}-x^{k-1}$ and integrating over (0,1), we obtain
\[\int\limits_0^1 {\left( {{x^{r - 1}} - {x^{k - 1}}} \right)\frac{{{\rm Li}{_m}\left( x \right)}}{{1 - x}}dx}  = \left( {k - r} \right)\sum\limits_{n = 1}^\infty  {\frac{{{\zeta _n}\left( m \right)}}{{\left( {n + r} \right)\left( {n + k} \right)}}} \;\;\;\left( {0 \le r < k,\;r,k \in \Z} \right),\tag{2.5}\]
\[\int\limits_0^1 {\left( {{x^{k - 1}} - {x^{r - 1}}} \right)\frac{{{\rm Li}{_m}\left( { - x} \right)}}{{1 - x}}} dx = \left( {k - r} \right)\sum\limits_{n = 1}^\infty  {\frac{{{L_n}\left( m \right)}}{{\left( {n + r} \right)\left( {n + k} \right)}}} \:\:\:\left( {0 \le r < k,\:r,k \in \Z} \right).\tag{2.6}\]
We now evaluate the integral on the left side of (2.5) and (2.6). Noting that
\[\int\limits_0^1 {\left( {{x^{r - 1}} - {x^{k - 1}}} \right)\frac{{{\rm Li}{_m}\left( x \right)}}{{1 - x}}dx}  = \sum\limits_{i = 1}^{k - r} {\int\limits_0^1 {{x^{r + i - 2}}{\rm Li}{_m}\left( x \right)dx}}, \tag{2.7}\]
\[\int\limits_0^1 {\left( {{x^{k - 1}} - {x^{r - 1}}} \right)\frac{{{\rm Li}{_m}\left( { - x} \right)}}{{1 - x}}} dx = \sum\limits_{i = 1}^{k - r} {{{\left( { - 1} \right)}^{r + i}}\int\limits_0^{ - 1} {{x^{r + i - 2}}{\rm Li}{_m}\left( x \right)} } dx.\tag{2.8}\]
Using integration by parts we have
\[\int\limits_0^x {{t^{n - 1}}{\rm Li}{_q}\left( t \right)dt}  = \sum\limits_{i = 1}^{q - 1} {{{\left( { - 1} \right)}^{i - 1}}\frac{{{x^n}}}{{{n^i}}}{\rm Li}{_{q + 1 - i}}\left( x \right)}  + \frac{{{{\left( { - 1} \right)}^q}}}{{{n^q}}}\ln \left( {1 - x} \right)\left( {{x^n} - 1} \right) - \frac{{{{\left( { - 1} \right)}^q}}}{{{n^q}}}\left( {\sum\limits_{k = 1}^n {\frac{{{x^k}}}{k}} } \right).\tag{2.9}\]
Taking $r=0$ in (2.5)-(2.8), substituting (2.9) into (2.5)-(2.8), we can obtain (2.1) and (2.2).\hfill$\square$ \\
If taking $r\geq 1$ in (2.5)-(2.8), using (2.9), we have
\[\sum\limits_{n = 1}^\infty  {\frac{{{\zeta _n}\left( m \right)}}{{\left( {n + r} \right)\left( {n + k} \right)}}}  = \frac{1}{{k - r}}\left\{ \begin{array}{l}
 \sum\limits_{j = 1}^{m - 1} {{{\left( { - 1} \right)}^{j - 1}}\zeta \left( {m + 1 - j} \right)\left( {{\zeta _{k - 1}}\left( j \right) - {\zeta _{r - 1}}\left( j \right)} \right)}  \\
  + {\left( { - 1} \right)^{m - 1}}\left( {\sum\limits_{i = 1}^{k - 1} {\frac{{{H_i}}}{{{i^m}}}}  - \sum\limits_{i = 1}^{r - 1} {\frac{{{H_i}}}{{{i^m}}}} } \right) \\
 \end{array} \right\},\tag{2.10}\]
\[\sum\limits_{n = 1}^\infty  {\frac{{{L_n}\left( m \right)}}{{\left( {n + r} \right)\left( {n + k} \right)}}}  = \frac{1}{{k - r}}\left\{ \begin{array}{l}
 \sum\limits_{j = 1}^{m - 1} {{{\left( { - 1} \right)}^{j - 1}}\bar \zeta \left( {m + 1 - j} \right)\left( {{\zeta _{k - 1}}\left( j \right) - {\zeta _{r - 1}}\left( j \right)} \right)}  \\
  + {\left( { - 1} \right)^{m - 1}}\ln 2\left( {{\zeta _{k - 1}}\left( m \right) - {\zeta _{r - 1}}\left( m \right) + {L_{k - 1}}\left( m \right) - {L_{r - 1}}\left( m \right)} \right) \\
  + {\left( { - 1} \right)^m}\sum\limits_{i = r}^{k - 1} {\frac{{{{\left( { - 1} \right)}^{i - 1}}}}{{{i^m}}}{L_i}\left( 1 \right)}  \\
 \end{array} \right\}.\tag{2.11}\]
Letting $m=1$ in (2.1) and (2.2), we conclude that
\begin{align*}
\sum\limits_{n = 1}^\infty  {\frac{{{H_n}}}{{n\left( {n + k} \right)}}}  = \frac{1}{k}\left( {\frac{1}{2}H_k^2 + \frac{1}{2}{\zeta _k}\left( 2 \right) + \zeta \left( 2 \right) - \frac{{{H_k}}}{k}} \right),\tag{2.12}
\end{align*}
\[\sum\limits_{n = 1}^\infty  {\frac{{{L_n}\left( 1 \right)}}{{n\left( {n + k} \right)}}}  = \frac{1}{k}\left\{ \begin{array}{l}
 \bar \zeta \left( 2 \right) + \ln 2\left( {{H_k} + {L_k}\left( 1 \right)} \right) - \ln 2\frac{{1 + {{\left( { - 1} \right)}^{k - 1}}}}{k} \\
  - \frac{1}{2}\left( {L_k^2\left( 1 \right) + {\zeta _k}\left( 2 \right)} \right) + \frac{{{L_k}\left( 1 \right)}}{k}{\left( { - 1} \right)^{k - 1}} \\
 \end{array} \right\}.\tag{2.13}\]
\begin{thm}\label{thm 2.1}
For integers $n\geq1, k\geq0$, we have
\[\int\limits_0^1 {{t^{n - 1}}{{\ln }^k}\left( {1 - t} \right)} dt = {\left( { - 1} \right)^k}\frac{{{Y_k}\left( n \right)}}{n},\tag{2.14}\]
where ${Y_k}\left( n \right) = {Y_k}\left( {{\zeta _n}\left( 1 \right),1!{\zeta _n}\left( 2 \right),2!{\zeta _n}\left( 3 \right), \cdots ,\left( {r - 1} \right)!{\zeta _n}\left( r \right), \cdots } \right)$, ${Y_k}\left( {{x_1},{x_2}, \cdots } \right)$ stands for the complete exponential Bell polynomial is defined by (see [15])
\[\exp \left( {\sum\limits_{m \ge 1}^{} {{x_m}\frac{{{t^m}}}{{m!}}} } \right) = 1 + \sum\limits_{k \ge 1}^{} {{Y_k}\left( {{x_1},{x_2}, \cdots } \right)\frac{{{t^k}}}{{k!}}}.\tag{2.15}\]
\end{thm}
\pf Using the definition of the complete exponential Bell polynomial, it is easily shown that
\[1 + \sum\limits_{k \ge 1}^{} {{Y_k}\left( n \right)\frac{{{t^k}}}{{k!}}}  = \frac{1}{{\left( {1 - t} \right)\left( {1 - \frac{t}{2}} \right) \cdots \left( {1 - \frac{t}{n}} \right)}}\]
and
\[\frac{1}{{\left( {1 - t} \right)\left( {1 - \frac{t}{2}} \right) \cdots \left( {1 - \frac{t}{n}} \right)}} = 1 + \sum\limits_{m = 1}^n {\frac{t}{{m\left( {1 - t} \right)\left( {1 - \frac{t}{2}} \right) \cdots \left( {1 - \frac{t}{m}} \right)}}} \]
Therefore we obtain
\[{Y_k}\left( n \right) = k\sum\limits_{m = 1}^n {\frac{{{Y_{k - 1}}\left( m \right)}}{m}} ,\;{Y_0}\left( n \right) = 1.\tag{2.16}\]
It is easily shown using integration by parts that
$$\int_0^x {{t^{n - 1}}\ln \left( {1 - t} \right)} dt = \frac {1}{n}\left\{{x^n}\ln \left( {1 - x} \right) - \sum\limits_{j = 1}^n {\frac{{{x^j}}}{j}}  - \ln \left( {1 - x} \right)\right\},  -1\leq x<1.$$
Letting $x\rightarrow 1^-$, we have
\[\mathop {\lim }\limits_{x \to {1^ - }} \int_0^x {{t^{n - 1}}\ln \left( {1 - t} \right)} dt =  - \frac{H_n}{n}=-\frac {Y_1(n)}{n}.\tag{2.17}\]
Using integration by parts and (2.17) , we can find that
\[\int_0^x {{t^{n - 1}}{{\ln }^2}\left( {1 - t} \right)} dt = \frac{1}{n} {{(x^n-1)}{{\ln }^2}\left( {1 - x} \right) } - \frac{2}{n}\sum\limits_{k = 1}^n {\frac{1}{k}} \left\{ {{x^k}\ln \left( {1 - x} \right) - \sum\limits_{j = 1}^k {\frac{{{x^j}}}{j}}  - \ln \left( {1 - x} \right)} \right\}\tag{2.18}\]
Letting $x$ tend to $ 1^-$ in (2.18), we deduce that
\[\mathop {\lim }\limits_{x \to {1^ - }} \int_0^x {{t^{n - 1}}{{\ln }^2}\left( {1 - t} \right)} dt = \frac{2}{n}\sum\limits_{k = 1}^n {\frac{1}{k}\sum\limits_{j = 1}^k {\frac{1}{j}} }  = \frac{{{Y_2}\left( n \right)}}{n}.\tag{2.19}\]
Further, it is easily verify that
\[\mathop {\lim }\limits_{x \to {1^ - }} \int_0^x {{t^{n - 1}}{{\ln }^m}\left( {1 - t} \right)} dt = m!\frac{(-1)^{m}}{n}\sum\limits_{{k_1} = 1}^n {\frac{1}{{{k_1}}}\sum\limits_{{k_2} = 1}^{{k_1}} {\frac{1}{{{k_2}}} \cdots } } \sum\limits_{{k_m} = 1}^{{k_{m - 1}}} {\frac{1}{{{k_m}}}}  = {\left( { - 1} \right)^m}\frac{{{Y_m}\left( n \right)}}{n}, 1\leq m\in \Z.\tag{2.20}\]
We complete the proof of (2.14).\hfill$\square$ \\
 From the definition of the complete exponential Bell polynomial, we have
$${Y_1}\left( n \right) = {H_n},{Y_2}\left( n \right) = H_n^2 + {\zeta _n}\left( 2 \right),{Y_3}\left( n \right) =  H_n^3+ 3{H_n}{\zeta _n}\left( 2 \right)+ 2{\zeta _n}\left( 3 \right),$$
\[{Y_4}\left( n \right) = H_n^4 + 8{H_n}{\zeta _n}\left( 3 \right) + 6H_n^2{\zeta _n}\left( 2 \right) + 3\zeta _n^2\left( 2 \right) + 6{\zeta _n}\left( 4 \right).\]
\begin{thm}\label{thm 2.2}\ \ For integer $p\geq 2,k\geq 1$, we have
\[\left( {p - 1} \right)!\sum\limits_{n = p - 1}^\infty  {\frac{{S\left( {n + 1,p} \right)}}{{n!n\left( {n + k} \right)}}}  = \frac{1}{k}\left\{ {\left( {p - 1} \right)!\zeta(p) + \frac{{{Y_p}\left( k \right)}}{p} - \frac{{{Y_{p - 1}}\left( k \right)}}{k}} \right\}.\tag{2.21}\]
where ${S\left( {n,k} \right)}$ are the Stirling number of the first kind defined by
\[n!\left( {1 + x} \right)\left( {1 + \frac{x}{2}} \right) \cdots \left( {1 + \frac{x}{n}} \right) = \sum\limits_{k = 0}^n {S\left( {n + 1,k + 1} \right){x^k}} .\]
From the definition ${S\left( {n,k} \right)}$, we can rewrite it as
\begin{align*}
\sum\limits_{k = 0}^n {S\left( {n + 1,k + 1} \right){x^k}}  &= n!\exp \left\{ {\sum\limits_{j = 1}^n {\ln \left( {1 + \frac{x}{j}} \right)} } \right\}\\
& = n!\exp \left\{ {\sum\limits_{j = 1}^n {\sum\limits_{k = 1}^\infty  {{{\left( { - 1} \right)}^{k - 1}}\frac{{{x^k}}}{{k{j^k}}}} } } \right\}\\
& = n!\exp \left\{ {\sum\limits_{k = 1}^\infty  {{{\left( { - 1} \right)}^{k - 1}}\frac{{{\zeta _n}\left( k \right){x^k}}}{k}} } \right\}.
\end{align*}
Therefore, we know that ${S\left( {n,k} \right)}$ is a rational linear combination of products of harmonic numbers.
The following identities is easily derived
\begin{align*}
& S\left( {n,1} \right) = \left( {n - 1} \right)!,S\left( {n,2} \right) = \left( {n - 1} \right)!{H_{n - 1}},S\left( {n,3} \right) = \frac{{\left( {n - 1} \right)!}}{2}\left[ {H_{n - 1}^2 - {\zeta _{n - 1}}\left( 2 \right)} \right],\\
&S\left( {n,4} \right) = \frac{{\left( {n - 1} \right)!}}{6}\left[ {H_{n - 1}^3 - 3{H_{n - 1}}{\zeta _{n - 1}}\left( 2 \right) + 2{\zeta _{n - 1}}\left( 3 \right)} \right], \\
&S\left( {n,5} \right) = \frac{{\left( {n - 1} \right)!}}{{24}}\left[ {H_{n - 1}^4 - 6{\zeta _{n - 1}}\left( 4 \right) - 6H_{n - 1}^2{\zeta _{n - 1}}\left( 2 \right) + 3\zeta _{n - 1}^2\left( 2 \right) + 8H_{n - 1}^{}{\zeta _{n - 1}}\left( 3 \right)} \right].
\end{align*}
\end{thm}
\pf From [15], we have \[{\ln ^p}\left( {1 - x} \right) = {\left( { - 1} \right)^p}p!\sum\limits_{n = p}^\infty  {S\left( {n,p} \right)} \frac{{{x^n}}}{{n!}}\:\:,1 \leq p\in \Z, -1\leq x<1.\tag{2.22}\]
Differentiating this equality, we obtain
\[\frac{{{{\ln }^{p - 1}}\left( {1 - x} \right)}}{{1 - x}} = {\left( { - 1} \right)^{p - 1}}\left( {p - 1} \right)!\sum\limits_{n = p - 1}^\infty  {S\left( {n + 1,p} \right)} \frac{{{x^n}}}{{n!}}\:\:,p \ge 2.\tag{2.23}\]
Let $k,r$ be integers with $k>r\geq 0$, then multiplying (2.23) by $x^{r-1}-x^{k-1}$ and integrating over (0,1), we get
\[\int\limits_0^1 {\frac{{{{\ln }^{p - 1}}\left( {1 - x} \right)}}{{1 - x}}\left( {{x^{r - 1}} - {x^{k - 1}}} \right)dx}  = {\left( { - 1} \right)^{p - 1}}\left( {p - 1} \right)!\sum\limits_{n = p - 1}^\infty  {\frac{{\left( {k - r} \right)S\left( {n + 1,p} \right)}}{{n!\left( {n + r} \right)\left( {n + k} \right)}}} \:\:,p \ge 2.\tag{2.24}\]
Noting that
\[\int\limits_0^1 {\frac{{{{\ln }^{p - 1}}\left( {1 - x} \right)}}{{1 - x}}\left( {{x^{r - 1}} - {x^{k - 1}}} \right)dx}  = \sum\limits_{i = 1}^{k - r} {\int\limits_0^1 {{x^{r + i - 2}}{{\ln }^{p - 1}}\left( {1 - x} \right)dx} } \tag{2.25}\]
and
\[\int\limits_0^1 {\frac{{{{\ln }^{p - 1}}\left( {1 - x} \right)}}{x}dx}  = {\left( { - 1} \right)^{p - 1}}\left( {p - 1} \right)!\zeta \left( p \right),\;2 \le p \in Z.\tag{2.26}\]
Taking $r=0, k\geq 1$ in (2.24), using (2.14)(2.16)(2.25) and (2.26), we can obtain (2.21).\hfill$\square$ \\ If taking $r\geq 1$ in (2.24), then we have
\[\left( {p - 1} \right)!\left( {k - r} \right)\sum\limits_{n = p - 1}^\infty  {\frac{{S\left( {n + 1,p} \right)}}{{n!\left( {n + r} \right)\left( {n + k} \right)}}}  = \frac{1}{p}\left\{ {{Y_p}\left( {k - 1} \right) - {Y_p}\left( {r - 1} \right)} \right\}.\tag{2.27}\]
Letting $p=3,4$ in (2.8), we can give the following Corollary.
\begin{cor}
If $1\leq k\in \Z$, then we have
\[\sum\limits_{n = 1}^\infty  {\frac{{H_n^2 - {\zeta _n}\left( 2 \right)}}{{n\left( {n + k} \right)}}}  = \frac{1}{k}\left\{ {2\zeta \left( 3 \right) + \frac{{H_k^3 + 3{H_k}{\zeta _k}\left( 2 \right) + 2{\zeta _k}\left( 3 \right)}}{3} - \frac{{H_k^2 + {\zeta _k}\left( 2 \right)}}{k}} \right\},\tag{2.28}\]
\[\sum\limits_{n = 1}^\infty  {\frac{{H_n^3 - 3{H_n}{\zeta _n}\left( 2 \right) + 2{\zeta _n}\left( 3 \right)}}{{n\left( {n + k} \right)}}}  = \frac{1}{k}\left\{ \begin{array}{l}
 \frac{{H_k^4 + 8{H_k}{\zeta _k}\left( 3 \right) + 6H_k^2{\zeta _k}\left( 2 \right) + 3\zeta _k^2\left( 2 \right) + 6{\zeta _k}\left( 4 \right)}}{4} \\
  - \frac{{H_k^3 + 3{H_k}{\zeta _k}\left( 2 \right) + 2{\zeta _k}\left( 3 \right)}}{k} + 6\zeta \left( 4 \right) \\
 \end{array} \right\}.\tag{2.29}\]
\end{cor}
From (2.1), taking $m=2$, we have
\[\sum\limits_{n = 1}^\infty  {\frac{{{\zeta _n}\left( 2 \right)}}{{n\left( {n + k} \right)}}}  = \frac{1}{k}\left\{ {\zeta \left( 3 \right) + \zeta \left( 2 \right){H_{k - 1}} - \sum\limits_{i = 1}^{k - 1} {\frac{{{H_i}}}{{{i^2}}}} } \right\}.\tag{2.30}\]
Substituting (2.30) into (2.28), we obtain (1.2).
In the same manner, we obtain the following Theorem:
\begin{thm} For integer $m,k>0$, then
\begin{align*}
&\sum\limits_{n = 1}^\infty  {\frac{{{\zeta _n}\left( m \right)}}{{n + k}}{{\left( { - 1} \right)}^{n - 1}}}  \\
& ={\left( { - 1} \right)^k}\left( {\sum\limits_{n = 1}^\infty  {\frac{{{\zeta _n}\left( m \right)}}{n}{{\left( { - 1} \right)}^{n - 1}} - \bar \zeta \left( {m + 1} \right)} } \right) + {\left( { - 1} \right)^{k + m + 1}}\ln 2\left( {{\zeta _{k - 1}}\left( m \right) + {L_{k - 1}}\left( m \right)} \right)\\
& \quad + {\left( { - 1} \right)^k}\sum\limits_{j = 1}^{m - 1} {{{\left( { - 1} \right)}^{j - 1}}\bar \zeta \left( {m + 1 - j} \right){L_{k - 1}}\left( j \right)}  + {\left( { - 1} \right)^{m + k}}\sum\limits_{i = 1}^{k - 1} {\frac{{{L_i}\left( 1 \right)}}{{{i^m}}}}, \tag{2.31}\\
&\sum\limits_{n = 1}^\infty  {\frac{{{L_n}\left( m \right)}}{{n + k}}{{\left( { - 1} \right)}^{n - 1}}} \\
& ={\left( { - 1} \right)^{k - 1}}\left( {\zeta \left( {m + 1} \right) - \sum\limits_{n = 1}^\infty  {\frac{{{L_n}\left( m \right)}}{n}{{\left( { - 1} \right)}^{n - 1}}}  + {{\left( { - 1} \right)}^m}\sum\limits_{i = 1}^{k - 1} {\frac{{{H_i}}}{{{i^m}}}{{\left( { - 1} \right)}^{i - 1}}} } \right)\\
& \quad + {\left( { - 1} \right)^k}\sum\limits_{j = 1}^{m - 1} {{{\left( { - 1} \right)}^{j - 1}}\zeta \left( {m + 1 - j} \right){L_{k - 1}}\left( j \right)}  \tag{2.32}
\end{align*}
\end{thm}
From (2.31) and (2.32), we obtain
\[\sum\limits_{n = 1}^\infty  {\frac{{{H_n}}}{{n + k}}{{\left( { - 1} \right)}^{n - 1}}}  = {\left( { - 1} \right)^{k - 1}}\left( {\frac{1}{2}{{\ln }^2}2 - \ln 2\left( {{H_{k - 1}} + {L_{k - 1}}\left( 1 \right)} \right) + \sum\limits_{i = 1}^{k - 1} {\frac{{{L_i}\left( 1 \right)}}{i}} } \right), \tag{2.33}\]
\[\sum\limits_{n = 1}^\infty  {\frac{{{L_n}\left( 1 \right)}}{{n + k}}{{\left( { - 1} \right)}^{n - 1}}}  = {\left( { - 1} \right)^{k - 1}}\left( {\frac{{\zeta \left( 2 \right) - {{\ln }^2}2}}{2} - \sum\limits_{i = 1}^{k - 1} {\frac{{{H_i}}}{i}{{\left( { - 1} \right)}^{i - 1}}} } \right). \tag{2.34}\]
\begin{thm}
For $1 \le {l_1},{l_2},m \in \Z$ and $x,y,z \in \left[ { - 1,1} \right)$, we have the following relation
\begin{align*}
&\sum\limits_{n = 1}^\infty  {\frac{{{\zeta _n}\left( {{l_1},x} \right){\zeta _n}\left( {{l_2},y} \right)}}{{{n^m}}}{z^n}}  + \sum\limits_{n = 1}^\infty  {\frac{{{\zeta _n}\left( {{l_1},x} \right){\zeta _n}\left( {m,z} \right)}}{{{n^{{l_2}}}}}{y^n}}  + \sum\limits_{n = 1}^\infty  {\frac{{{\zeta _n}\left( {{l_2},y} \right){\zeta _n}\left( {m,z} \right)}}{{{n^{{l_1}}}}}{x^n}} \\
& =\sum\limits_{n = 1}^\infty  {\frac{{{\zeta _n}\left( {m,z} \right)}}{{{n^{{l_1} + {l_2}}}}}{{\left( {xy} \right)}^n}}  + \sum\limits_{n = 1}^\infty  {\frac{{{\zeta _n}\left( {{l_1},x} \right)}}{{{n^{m + {l_2}}}}}{{\left( {yz} \right)}^n}}  + \sum\limits_{n = 1}^\infty  {\frac{{{\zeta _n}\left( {{l_2},y} \right)}}{{{n^{{l_1} + m}}}}{{\left( {xz} \right)}^n}}  \\
&\quad  +{\rm Li}{_m}\left( z \right){\rm Li}{_{{l_1}}}\left( x \right){\rm Li}{_{{l_2}}}\left( y \right) - {\rm Li}{_{{l_1} + {l_2} + m}}\left( {xyz} \right) \tag{2.35}
\end{align*}
where the partial sum ${\zeta _n}\left( {l,x} \right)$ of polylogarithm function is defined by ${\zeta _n}\left( {l,x} \right) := \sum\limits_{k = 1}^n {\frac{{{x^k}}}{k^l}}$.
\end{thm}
\pf  We construct the function $F\left( {x,y,z} \right) := \sum\limits_{n = 1}^\infty  {\left\{ {{\zeta _n}\left( {{l_1},x} \right){\zeta _n}\left( {{l_2},y} \right) - {\zeta _n}\left( {{l_1} + {l_2},xy} \right)} \right\}{z^{n - 1}}} $, $z\in (-1,1)$. By the definition of ${\zeta _n}\left( {l,x} \right)$, we have
\[F\left( {x,y,z} \right) = zF\left( {x,y,z} \right) + \sum\limits_{n = 1}^\infty  {\left\{ {\frac{{{\zeta _n}\left( {{l_1},x} \right)}}{{{{\left( {n + 1} \right)}^{{l_2}}}}}{y^{n+1}} + \frac{{{\zeta _n}\left( {{l_2},y} \right)}}{{{{\left( {n + 1} \right)}^{{l_1}}}}}{x^{n+1}}} \right\}{z^n}} .\tag{2.36}\]
Moving $zF(x,y,z)$ from right to left and then multiplying  $(1-z)^{-1}$  to the equation (2.36) and integrating over the interval $(0,z)$, we obtain
\begin{align*}
&\sum\limits_{n = 1}^\infty  {\frac{{{\zeta _n}\left( {{l_1},x} \right){\zeta _n}\left( {{l_2},y} \right) - {\zeta _n}\left( {{l_1} + {l_2},xy} \right)}}{n}{z^{n}}} \\ &= \sum\limits_{n = 1}^\infty  {\left\{ {\frac{{{\zeta _n}\left( {{l_1},x} \right)}}{{{{\left( {n + 1} \right)}^{{l_2}}}}}{y^{n+1}} + \frac{{{\zeta _n}\left( {{l_2},y} \right)}}{{{{\left( {n + 1} \right)}^{{l_1}}}}}{x^{n+1}}} \right\}\left\{ {{\rm Li}{_1}\left( z \right) - {\zeta _n}\left( {1,z} \right)} \right\}}. \tag{2.37}
\end{align*}
Furthermore, using integration and the following formula
\[\sum\limits_{n = 1}^\infty  {\left\{ {\frac{{{\zeta _n}\left( {{l_1},x} \right)}}{{{{\left( {n + 1} \right)}^{{l_2}}}}}{y^{n + 1}} + \frac{{{\zeta _n}\left( {{l_2},y} \right)}}{{{{\left( {n + 1} \right)}^{{l_1}}}}}{x^{n + 1}}} \right\}}  = {\rm Li}{_{{l_1}}}\left( x \right){\rm Li}{_{{l_2}}}\left( y \right) - {\rm Li}{_{{l_1} + {l_2}}}\left( {xy} \right),\]
we deduce (2.35) to complete the proof of Theorem 2.6. \hfill$\square$ \\
Letting $(x,y,z)=(-1,-1,1),\ (l_1,l_2,m)=(1,1,m)$ and $(x,y,z)=(-1,-1,-1),\ (l_1,l_2,m)=(1,1,m)$ in (2.35) gives the following
\begin{cor} For integer $m>1$ we have that
\begin{align*}
&\sum\limits_{n = 1}^\infty  {\frac{{L_n^2\left( 1 \right)}}{{{n^m}}}}  + 2\sum\limits_{n = 1}^\infty  {\frac{{{L_n}\left( 1 \right){\zeta _n}\left( m \right)}}{n}{{\left( { - 1} \right)}^{n - 1}}}\\ &=\sum\limits_{n = 1}^\infty  {\frac{{{\zeta _n}\left( m \right)}}{{{n^2}}}}  + 2\sum\limits_{n = 1}^\infty  {\frac{{{L_n}\left( 1 \right)}}{{{n^{m + 1}}}}{{\left( { - 1} \right)}^{n - 1}}}  + {\ln ^2}2\zeta \left( m \right) - \zeta \left( {m + 2} \right), \ (m>1)\tag{2.38}\\
&\sum\limits_{n = 1}^\infty  {\frac{{L_n^2\left( 1 \right)}}{{{n^m}}}{{\left( { - 1} \right)}^{n - 1}}}  + 2\sum\limits_{n = 1}^\infty  {\frac{{{L_n}\left( 1 \right){L_n}\left( m \right)}}{n}{{\left( { - 1} \right)}^{n - 1}}}\\ &=\sum\limits_{n = 1}^\infty  {\frac{{{L_n}\left( m \right)}}{{{n^2}}}}  + 2\sum\limits_{n = 1}^\infty  {\frac{{{L_n}\left( 1 \right)}}{{{n^{m + 1}}}}}  + {\ln ^2}2\bar \zeta \left( m \right) - \bar \zeta \left( {m + 2} \right), \ (m>0). \tag{2.39}
\end{align*}
\end{cor}
In fact, multiplying (2.35) by $(1-z)^{-1}$ and integrating over $(0,z)\ (z\in (-1,1))$, we can obtain the following Corollary.
\begin{cor} For positive integers $l_1>0,l_2>0,m>1$ and $x,y,z\in [-1,1)$, then we have
\begin{align*}
&\sum\limits_{n = 1}^\infty  { {\frac{{{\zeta _n}\left( {{l_1},x} \right){\zeta _n}\left( {{l_2},x} \right){\zeta _n}\left( {1,z} \right)}}{{{n^m}}} + \sum\limits_{n = 1}^\infty\frac{{{\zeta _n}\left( {{l_1},x} \right)\left( {\sum\limits_{k = 1}^n {\frac{{{\zeta _k}\left( {1,z} \right)}}{{{k^m}}}} } \right)}}{{{n^{{l_2}}}}}{y^n} + \sum\limits_{n = 1}^\infty \frac{{{\zeta _n}\left( {{l_2},y} \right)\left( {\sum\limits_{k = 1}^n {\frac{{{\zeta _k}\left( {1,z} \right)}}{{{k^m}}}} } \right)}}{{{n^{{l_1}}}}}{x^n}} } \\
& = \sum\limits_{n = 1}^\infty  { {\frac{{\left( {\sum\limits_{k = 1}^n {\frac{{{\zeta _k}\left( {1,z} \right)}}{{{k^m}}}} } \right)}}{{{n^{{l_1} + {l_2}}}}}{x^n}{y^n} + \sum\limits_{n = 1}^\infty\frac{{{\zeta _n}\left( {{l_1},x} \right){\zeta _n}\left( {1,z} \right)}}{{{n^{m + {l_2}}}}}{y^n} + \sum\limits_{n = 1}^\infty\frac{{{\zeta _n}\left( {{l_2},y} \right){\zeta _n}\left( {1,z} \right)}}{{{n^{m + {l_1}}}}}{x^n}} } \\
&\quad + {\rm{L}}{{\rm{i}}_{{l_1}}}\left( x \right){\rm{L}}{{\rm{i}}_{{l_2}}}\left( y \right)\left( {\sum\limits_{n = 1}^\infty  {\frac{{{\zeta _n}\left( {1,z} \right)}}{{{n^m}}}} } \right) - \left( {\sum\limits_{n = 1}^\infty  {\frac{{{\zeta _n}\left( {1,z} \right)}}{{{n^{m + {l_1} + {l_2}}}}}{x^n}{y^n}} } \right). \tag{2.40}
\end{align*}
\end{cor}
Letting $x,y,z\rightarrow 1$ in (2.40), we arrive at the conclusion that, for integers $m,l_1,l_2>1$
\begin{align*}
&\sum\limits_{n = 1}^\infty  { {\frac{{{H_n}{\zeta _n}\left( {{l_1}} \right){\zeta _n}\left( {{l_2}} \right)}}{{{n^m}}} +\sum\limits_{n = 1}^\infty \frac{{{\zeta _n}\left( {{l_1}} \right)\left( {\sum\limits_{k = 1}^n {\frac{{{H_k}}}{{{k^m}}}} } \right)}}{{{n^{{l_2}}}}} + \sum\limits_{n = 1}^\infty\frac{{{\zeta _n}\left( {{l_2},y} \right)\left( {\sum\limits_{k = 1}^n {\frac{{{H_k}}}{{{k^m}}}} } \right)}}{{{n^{{l_1}}}}}} } \\
&= \sum\limits_{n = 1}^\infty  { {\frac{{\left( {\sum\limits_{k = 1}^n {\frac{{{H_k}}}{{{k^m}}}} } \right)}}{{{n^{{l_1} + {l_2}}}}} +\sum\limits_{n = 1}^\infty \frac{{{H_n}{\zeta _n}\left( {{l_1}} \right)}}{{{n^{m + {l_2}}}}} +\sum\limits_{n = 1}^\infty \frac{{{H_n}{\zeta _n}\left( {{l_2}} \right)}}{{{n^{m + {l_1}}}}}} } + \zeta \left( {{l_1}} \right)\zeta \left( {{l_2}} \right)\left( {\sum\limits_{n = 1}^\infty  {\frac{{{H_n}}}{{{n^m}}}} } \right) - \left( {\sum\limits_{n = 1}^\infty  {\frac{{{H_n}}}{{{n^{m + {l_1} + {l_2}}}}}} } \right) .\tag{2.41}
\end{align*}
From [13], we derive the following identity
\[\sum\limits_{n = 1}^\infty  {\frac{1}{{{n^{p + 1}}}}\left( {\sum\limits_{k = 1}^n {\frac{{{H_k}}}{{{k^m}}}} } \right)}  = \sum\limits_{n = 1}^\infty  {\frac{{{H_n}}}{{{n^{m + p + 1}}}}}  + \zeta \left( {p + 1} \right)\sum\limits_{n = 1}^\infty  {\frac{{{H_n}}}{{{n^m}}}}  - \sum\limits_{n = 1}^\infty  {\frac{{{H_n}{\zeta _n}\left( {p + 1} \right)}}{{{n^m}}}}.\tag{2.42}\]
Therefore, we may rewrite (2.41) as
\begin{align*}
&\sum\limits_{n = 1}^\infty  {{\frac{{{H_n}{\zeta _n}\left( {{l_1}} \right){\zeta _n}\left( {{l_2}} \right)}}{{{n^m}}} + \sum\limits_{n = 1}^\infty\frac{{{\zeta _n}\left( {{l_1}} \right)\left( {\sum\limits_{k = 1}^n {\frac{{{H_k}}}{{{k^m}}}} } \right)}}{{{n^{{l_2}}}}} +\sum\limits_{n = 1}^\infty \frac{{{\zeta _n}\left( {{l_2},y} \right)\left( {\sum\limits_{k = 1}^n {\frac{{{H_k}}}{{{k^m}}}} } \right)}}{{{n^{{l_1}}}}}} } \\
& = \sum\limits_{n = 1}^\infty  {\left\{ {\frac{{{H_n}{\zeta _n}\left( {{l_1}} \right)}}{{{n^{m + {l_2}}}}} + \frac{{{H_n}{\zeta _n}\left( {{l_2}} \right)}}{{{n^{m + {l_1}}}}} - \frac{{{H_n}{\zeta _n}\left( {{l_1} + {l_2}} \right)}}{{{n^m}}}} \right\}} \\
&\quad+ \zeta \left( {{l_1} + {l_2}} \right)\left( {\sum\limits_{n = 1}^\infty  {\frac{{{H_n}}}{{{n^m}}}} } \right) + \zeta \left( {{l_1}} \right)\zeta \left( {{l_2}} \right)\left( {\sum\limits_{n = 1}^\infty  {\frac{{{H_n}}}{{{n^m}}}} } \right).\tag{2.43}
\end{align*}
Taking $l_1=l_2=m=2l+1$ ($l$ is a positive integer) in (2.43), we conclude that
\begin{align*}
&\sum\limits_{n = 1}^\infty  { {\frac{{{H_n}\zeta _n^2\left( {2l + 1} \right)}}{{{n^{2l + 1}}}} + 2\sum\limits_{n = 1}^\infty\frac{{{\zeta _n}\left( {2l + 1} \right)\left( {\sum\limits_{k = 1}^n {\frac{{{H_k}}}{{{k^{2l + 1}}}}} } \right)}}{{{n^{2l + 1}}}}} } \\
& = \sum\limits_{n = 1}^\infty  {\left\{ {2\frac{{{H_n}{\zeta _n}\left( {2l + 1} \right)}}{{{n^{4l + 2}}}} - \frac{{{H_n}{\zeta _n}\left( {4l + 2} \right)}}{{{n^{2l + 1}}}}} \right\}}  + \left( {\zeta \left( {4l + 2} \right) + {\zeta ^2}\left( {2l + 1} \right)} \right)\left( {\sum\limits_{n = 1}^\infty  {\frac{{{H_n}}}{{{n^{2l + 1}}}}} } \right).\tag{2.44}
\end{align*}
\section{Closed form of quadratic Euler sums}
In this section we evaluate some quadratic Euler sums involving harmonic numbers and alternating harmonic numbers.
\begin{thm} For integers $l>0,m>0,p>1$, we have
\begin{align*}
&{\left( { - 1} \right)^{m - 1}}\sum\limits_{n = 1}^\infty  {\frac{{{\zeta _n}\left( l \right)}}{{{n^{p + 1}}}}\left( {\sum\limits_{k = 1}^n {\frac{{{H_k}}}{{{k^m}}}} } \right)}  - {\left( { - 1} \right)^{p + l}}\sum\limits_{n = 1}^\infty  {\frac{{{\zeta _n}\left( m \right)}}{{{n^{p + 1}}}}\left( {\sum\limits_{k = 1}^n {\frac{{{H_k}}}{{{k^l}}}} } \right)} \\ &=\sum\limits_{i = 1}^{p - 1} {{{\left( { - 1} \right)}^{i - 1}}\left( {\sum\limits_{n = 1}^\infty  {\frac{{{\zeta _n}\left( m \right)}}{{{n^{i + 1}}}}} } \right)\left( {\sum\limits_{n = 1}^\infty  {\frac{{{\zeta _n}\left( l \right)}}{{{n^{p + 1 - i}}}}} } \right)}  + {\left( { - 1} \right)^{p - 1}}\zeta \left( {l + 1} \right)\left( {\sum\limits_{n = 1}^\infty  {\frac{{{\zeta _n}\left( m \right)}}{{{n^{p + 1}}}}} } \right)\\ &\quad
 + {\left( { - 1} \right)^{p - 1}}\sum\limits_{j = 1}^{l - 1} {{{\left( { - 1} \right)}^{j - 1}}\zeta \left( {l + 1 - j} \right)\left\{ {\sum\limits_{n = 1}^\infty  {\frac{{{\zeta _n}\left( m \right){\zeta _n}\left( j \right)}}{{{n^{p + 1}}}}}  - \sum\limits_{n = 1}^\infty  {\frac{{{\zeta _n}\left( m \right)}}{{{n^{p + j + 1}}}}} } \right\}}\\ &\quad
  - {\left( { - 1} \right)^{p + l}}\sum\limits_{n = 1}^\infty  {\frac{{{H_n}{\zeta _n}\left( m \right)}}{{{n^{p + l + 1}}}}}  + {\left( { - 1} \right)^{m - 1}}\sum\limits_{n = 1}^\infty  {\frac{{{H_n}{\zeta _n}\left( l \right)}}{{{n^{p + m + 1}}}}}  - \zeta \left( {m + 1} \right)\left( {\sum\limits_{n = 1}^\infty  {\frac{{{\zeta _n}\left( l \right)}}{{{n^{p + 1}}}}} } \right)\\ &\quad
   - \sum\limits_{j = 1}^{m - 1} {{{\left( { - 1} \right)}^{j - 1}}\zeta \left( {m + 1 - j} \right)\left\{ {\sum\limits_{n = 1}^\infty  {\frac{{{\zeta _n}\left( l \right){\zeta _n}\left( j \right)}}{{{n^{p + 1}}}}}  - \sum\limits_{n = 1}^\infty  {\frac{{{\zeta _n}\left( l \right)}}{{{n^{p + j + 1}}}}} } \right\}}.
 \tag{3.1}
\end{align*}
\end{thm}
\pf  Multiplying (2.1) by $\frac{\zeta_k(l)}{k^{p}}$ and summing with respect to $k$, we obtain
\[\sum\limits_{k = 1}^\infty  {\sum\limits_{n = 1}^\infty  {\frac{{{\zeta _k}\left( l \right){\zeta _n}\left( m \right)}}{{{k^p}n\left( {n + k} \right)}}} }  = \sum\limits_{k = 1}^\infty  {\frac{{{\zeta _k}\left( l \right)}}{{{k^p}}}} \sum\limits_{n = 1}^\infty  {\frac{{{\zeta _n}\left( m \right)}}{{n\left( {n + k} \right)}}}  = \sum\limits_{n = 1}^\infty  {\frac{{{\zeta _n}\left( m \right)}}{n}} \sum\limits_{k = 1}^\infty  {\frac{{{\zeta _k}\left( l \right)}}{{{k^p}\left( {n + k} \right)}}} .\]
Then by using (2.1) and the following partial fraction decomposition formula
\[\frac{1}{{{k^p}\left( {n + k} \right)}} = \sum\limits_{i = 1}^{p - 1} {\frac{{{{\left( { - 1} \right)}^{i - 1}}}}{{{n^i}}}\frac{1}{{{k^{p + 1 - i}}}}}  + \frac{{{{\left( { - 1} \right)}^{p - 1}}}}{{{n^{p - 1}}}}\frac{1}{{k\left( {n + k} \right)}},\]
we can obtain (3.1).\hfill$\square$ \\
Taking $(p,l)=(2l,m)$ in (3.1), we can find that
\begin{align*}
& \sum\limits_{n = 1}^\infty  {\frac{{{\zeta _n}\left( m \right)}}{{{n^{2l + 1}}}}\left( {\sum\limits_{k = 1}^n {\frac{{{H_k}}}{{{k^m}}}} } \right)}\\
 & = {\left( { - 1} \right)^{m - 1}}\sum\limits_{i = 1}^l {{{\left( { - 1} \right)}^{i - 1}}\left( {\sum\limits_{n = 1}^\infty  {\frac{{{\zeta _n}\left( m \right)}}{{{n^{i + 1}}}}} } \right)\left( {\sum\limits_{n = 1}^\infty  {\frac{{{\zeta _n}\left( m \right)}}{{{n^{2l + 1 - i}}}}} } \right)}  + \frac{{{{\left( { - 1} \right)}^{m + l - 1}}}}{2}{\left( {\sum\limits_{n = 1}^\infty  {\frac{{{\zeta _n}\left( m \right)}}{{{n^{l + 1}}}}} } \right)^2} \\
 & \quad- {\left( { - 1} \right)^{m - 1}}\zeta \left( {m + 1} \right)\left( {\sum\limits_{n = 1}^\infty  {\frac{{{\zeta _n}\left( m \right)}}{{{n^{2l + 1}}}}} } \right) + \sum\limits_{n = 1}^\infty  {\frac{{{H_n}{\zeta _n}\left( m \right)}}{{{n^{2l + m + 1}}}}}  \\
& \quad - {\left( { - 1} \right)^{m - 1}}\sum\limits_{j = 1}^{m - 1} {{{\left( { - 1} \right)}^{j - 1}}\zeta \left( {m + 1 - j} \right)\left\{ {\sum\limits_{n = 1}^\infty  {\frac{{{\zeta _n}\left( m \right){\zeta _n}\left( j \right)}}{{{n^{2l + 1}}}}}  - \sum\limits_{n = 1}^\infty  {\frac{{{\zeta _n}\left( m \right)}}{{{n^{2l + j + 1}}}}} } \right\}} .\tag{3.2}
 \end{align*}
Putting $m=2l+1$ in (3.2) and combining (2.44), we obtain
\begin{align*}
{S_{1{{\left( {2l + 1} \right)}^2},0,\left( {2l + 1} \right)}} = &\sum\limits_{n = 1}^\infty  {\frac{{{H_n}\zeta _n^2\left( {2l + 1} \right)}}{{{n^{2l + 1}}}}} \\
 = &2\zeta \left( {2l + 2} \right)\left( {\sum\limits_{n = 1}^\infty  {\frac{{{\zeta _n}\left( {2l + 1} \right)}}{{{n^{2l + 1}}}}} } \right)\\
&+ \left( {\zeta \left( {4l + 2} \right) + {\zeta ^2}\left( {2l + 1} \right)} \right)\left( {\sum\limits_{n = 1}^\infty  {\frac{{{H_n}}}{{{n^{2l + 1}}}}} } \right)\\
& - {\left( { - 1} \right)^l}{\left( {\sum\limits_{n = 1}^\infty  {\frac{{{\zeta _n}\left( {2l + 1} \right)}}{{{n^{l + 1}}}}} } \right)^2} - \sum\limits_{n = 1}^\infty  {\frac{{{H_n}{\zeta _n}\left( {4l + 2} \right)}}{{{n^{2l + 1}}}}} \\
& - 2\sum\limits_{i = 1}^l {{{\left( { - 1} \right)}^i}\left( {\sum\limits_{n = 1}^\infty  {\frac{{{\zeta _n}\left( {2l + 1} \right)}}{{{n^{i + 1}}}}} } \right)\left( {\sum\limits_{n = 1}^\infty  {\frac{{{\zeta _n}\left( {2l + 1} \right)}}{{{n^{2l + 1 - i}}}}} } \right)} \\
& + 2\sum\limits_{j = 1}^{2l} {{{\left( { - 1} \right)}^{j - 1}}\zeta \left( {2l + 2 - j} \right)} \left\{ {\sum\limits_{n = 1}^\infty  {\frac{{{\zeta _n}\left( {2l + 1} \right){\zeta _n}\left( j \right)}}{{{n^{2l + 1}}}}}  - \sum\limits_{n = 1}^\infty  {\frac{{{\zeta _n}\left( {2l + 1} \right)}}{{{n^{2l + j + 1}}}}} } \right\}. \tag{3.3}
\end{align*}
Letting $l=1, p=m-1$ in (3.1) and using the formula
\[\sum\limits_{n = 1}^\infty  {\frac{{{H_n}}}{{{n^m}}}} \left( {\sum\limits_{k = 1}^n {\frac{{{H_k}}}{{{k^m}}}} } \right) = \frac{1}{2}\left\{ {{{\left( {\sum\limits_{n = 1}^\infty  {\frac{{{H_n}}}{{{n^m}}}} } \right)}^2} + \sum\limits_{n = 1}^\infty  {\frac{{H_n^2}}{{{n^{2m}}}}} } \right\},\;2 \le m \in \Z,\]
we can get the following result
\begin{align*}
{S_{{1^2}{m},0,m}}&=\sum\limits_{n = 1}^\infty  {\frac{{H_n^2{\zeta _n}\left( m \right)}}{{{n^m}}}} \\& = 2{\left( { - 1} \right)^{m - 1}}\sum\limits_{i = 1}^{m - 2} {{{\left( { - 1} \right)}^{i - 1}}\left( {\sum\limits_{n = 1}^\infty  {\frac{{{\zeta _n}\left( m \right)}}{{{n^{i + 1}}}}} } \right)\left( {\sum\limits_{n = 1}^\infty  {\frac{{{H_n}}}{{{n^{m - i}}}}} } \right)} \\
&\quad - \zeta \left( 2 \right)\left( {{\zeta ^2}\left( m \right) + \zeta \left( {2m} \right)} \right) - 2{\left( { - 1} \right)^{m - 1}}\zeta \left( {m + 1} \right)\left( {\sum\limits_{n = 1}^\infty  {\frac{{{H_n}}}{{{n^m}}}} } \right)\\
&\quad - 2{\left( { - 1} \right)^{m - 1}}\sum\limits_{j = 1}^{m - 1} {{{\left( { - 1} \right)}^{j - 1}}\zeta \left( {m + 1 - j} \right)\left\{ {\sum\limits_{n = 1}^\infty  {\frac{{{H_n}{\zeta _n}\left( j \right)}}{{{n^m}}}}  - \sum\limits_{n = 1}^\infty  {\frac{{{H_n}}}{{{n^{m + j}}}}} } \right\}}\\
&\quad + 2\sum\limits_{n = 1}^\infty  {\frac{{{H_n}{\zeta _n}\left( m \right)}}{{{n^{m + 1}}}}}  + \sum\limits_{n = 1}^\infty  {\frac{{H_n^2}}{{{n^{2m}}}}}  - {\left( {\sum\limits_{n = 1}^\infty  {\frac{{{H_n}}}{{{n^m}}}} } \right)^2} - \sum\limits_{n = 1}^\infty  {\frac{{{\zeta _n}\left( 2 \right){\zeta _n}\left( m \right)}}{{{n^m}}}} .\tag{3.4}
 \end{align*}
In [13], Philippe Flajolet and Bruno Salvy gave the following conclusion:
If $p_1+p_2+q$ is even, and $p_1>1, p_2>1, q>1$, the quadratic sums
\[{S_{{p_1}{p_2},0,q}} = \sum\limits_{n = 1}^\infty  {\frac{{{\zeta _n}\left( {{p_1}} \right){\zeta _n}\left( {{p_2}} \right)}}{{{n^q}}}} \] are reducible to linear sums.
In [27],  we showed that all quadratic Euler sums of the form
\[{S_{{1}{m},0,p}}=\sum\limits_{n = 1}^\infty  {\frac{{{H_n}{\zeta _n}\left( m \right)}}{{{n^p}}}}\ \ \left( {m + p \le 8} \right)\]
are reducible to polynomials in zeta values and to linear sums. Hence, from (3.4), we know that the cubic sums ${S_{{1^2}{m},0,m}}$ are reducible to polynomials in zeta values and to linear sums when $m=2,3,4,5$. For example
 \[\begin{array}{l}
 {S_{{1^2}2,0,2}} = \sum\limits_{n = 1}^\infty  {\frac{{H_n^2{\zeta _n}\left( 2 \right)}}{{{n^2}}}}  = \frac{{41}}{{12}}\zeta \left( 6 \right) + 2{\zeta ^2}\left( 3 \right), \\
 {S_{{1^2}3,0,3}} = \sum\limits_{n = 1}^\infty  {\frac{{H_n^2{\zeta _n}\left( 3 \right)}}{{{n^3}}}}  = \frac{9}{2}\zeta \left( 3 \right)\zeta \left( 5 \right) + \frac{3}{2}\zeta \left( 2 \right){\zeta ^2}\left( 3 \right) - \frac{{443}}{{288}}\zeta \left( 8 \right) - \frac{{23}}{4}{S_{2,0,6}}, \\
 \end{array}\]
where ${S_{2,0,6}} = \sum\limits_{n = 1}^\infty  {\frac{{{\zeta _n}\left( 2 \right)}}{{{n^6}}}} $.\\
Noting that, from Theorem 4.2 in the reference [13], we deduce that
\[{S_{23,0,3}} = \sum\limits_{n = 1}^\infty  {\frac{{{\zeta _n}\left( 2 \right){\zeta _n}\left( 3 \right)}}{{{n^3}}}}  = \frac{{45}}{2}\zeta \left( 3 \right)\zeta \left( 5 \right) - \frac{{827}}{{48}}\zeta \left( 8 \right) - \frac{3}{2}\zeta \left( 2 \right){\zeta ^2}\left( 3 \right) - \frac{{23}}{4}{S_{2,0,6}}.\]
In the same way, we can obtain the following Theorems.
\begin{thm} For integers $p_1>1,p_2>1,m>0$, we have
\begin{align*}
&\frac{{\left( {{p_1} - 1} \right)!}}{{{p_2}}}\sum\limits_{n = {p_1} - 1}^\infty  {\frac{{S\left( {n + 1,{p_1}} \right){Y_{{p_2}}}\left( n \right)}}{{{n^{m + 1}}n!}}}  - {\left( { - 1} \right)^{m - 1}}\frac{{\left( {{p_2} - 1} \right)!}}{{{p_1}}}\sum\limits_{n = {p_2} - 1}^\infty  {\frac{{S\left( {n + 1,{p_2}} \right){Y_{{p_1}}}\left( n \right)}}{{{n^{m + 1}}n!}}} \\
& = \left( {{p_1} - 1} \right)!\sum\limits_{n = {p_1} - 1}^\infty  {\frac{{S\left( {n + 1,{p_1}} \right){Y_{{p_2} - 1}}\left( n \right)}}{{{n^{m + 2}}n!}}}  - {\left( { - 1} \right)^{m - 1}}\left( {{p_2} - 1} \right)!\sum\limits_{n = {p_2} - 1}^\infty  {\frac{{S\left( {n + 1,{p_2}} \right){Y_{{p_1} - 1}}\left( n \right)}}{{{n^{m + 2}}n!}}}\\
& \quad+ \left( {{p_1} - 1} \right)!\left( {{p_2} - 1} \right)!\sum\limits_{i = 1}^{m - 1} {{{\left( { - 1} \right)}^{i - 1}}\left( {\sum\limits_{n = {p_1} - 1}^\infty  {\frac{{S\left( {n + 1,{p_1}} \right)}}{{{n^{m + 1 - i}}n!}}} } \right)\left( {\sum\limits_{n = {p_2} - 1}^\infty  {\frac{{S\left( {n + 1,{p_2}} \right)}}{{{n^{i + 1}}n!}}} } \right)} \\&\quad  + {\left( { - 1} \right)^{m - 1}}\left( {{p_1} - 1} \right)!\left( {{p_2} - 1} \right)!\zeta \left( {{p_1}} \right)\left( {\sum\limits_{n = {p_2} - 1}^\infty  {\frac{{S\left( {n + 1,{p_2}} \right)}}{{{n^{m + 1}}n!}}} } \right)\\ &
\quad- \left( {{p_1} - 1} \right)!\left( {{p_2} - 1} \right)!\zeta \left( {{p_2}} \right)\left( {\sum\limits_{n = {p_1} - 1}^\infty  {\frac{{S\left( {n + 1,{p_1}} \right)}}{{{n^{m + 1}}n!}}} } \right) .\tag{3.5}
\end{align*}
\end{thm}
\pf Replacing $p$ by $p_2$ in (2.21), we get
\[\left( {{p_2} - 1} \right)!\sum\limits_{n = {p_2} - 1}^\infty  {\frac{{S\left( {n + 1,{p_2}} \right)}}{{n!n\left( {n + k} \right)}}}  = \frac{1}{k}\left\{ {\left( {{p_2} - 1} \right)!\zeta ({p_2}) + \frac{{{Y_{{p_2}}}\left( k \right)}}{{{p_2}}} - \frac{{{Y_{{p_2} - 1}}\left( k \right)}}{k}} \right\} .\tag{3.6}\]
 Multiplying (3.6) by $\left( {{p_1} - 1} \right)!\frac{{S\left( {k + 1,{p_1}} \right)}}{{k!{k^m}}}$ and summing with respect to $k$, we obtain
\[\begin{array}{l}
 \left( {{p_1} - 1} \right)!\left( {{p_2} - 1} \right)!\sum\limits_{k = {p_1} - 1}^\infty  {\sum\limits_{n = {p_2} - 1}^\infty  {\frac{{S\left( {k + 1,{p_1}} \right)S\left( {n + 1,{p_2}} \right)}}{{k!{k^m}n!n\left( {n + k} \right)}}} }  \\
  = \left( {{p_1} - 1} \right)!\sum\limits_{k = {p_1} - 1}^\infty  {\frac{{S\left( {k + 1,{p_1}} \right)}}{{k!{k^m}}}} \left( {{p_2} - 1} \right)!\sum\limits_{n = {p_2} - 1}^\infty  {\frac{{S\left( {n + 1,{p_2}} \right)}}{{n!n\left( {n + k} \right)}}}  \\
  = \left( {{p_2} - 1} \right)!\sum\limits_{n = {p_2} - 1}^\infty  {\frac{{S\left( {n + 1,{p_2}} \right)}}{{n!n}}} \left( {{p_1} - 1} \right)!\sum\limits_{k = {p_1} - 1}^\infty  {\frac{{S\left( {k + 1,{p_1}} \right)}}{{k!{k^m}\left( {n + k} \right)}}} . \\
 \end{array}\]
Then with the help of formula (2.21) we may easily deduce the result.\hfill$\square$
\begin{thm} For integer $m>0$, then we have
\begin{align*}
&\left( {\frac{1}{2} + {{\left( { - 1} \right)}^m}} \right)\sum\limits_{n = 1}^\infty  {\frac{{H_n^2}}{{{n^{m + 1}}}}}\\
& = \sum\limits_{j = 0}^{m - 2} {{{\left( { - 1} \right)}^j}} \zeta \left( {m - j} \right)\sum\limits_{n = 1}^\infty  {\frac{{{H_n}}}{{{n^{j + 2}}}}}  - \zeta \left( 2 \right)\zeta \left( {m + 1} \right) + \sum\limits_{n = 1}^\infty  {\frac{{H_n^{}}}{{{n^{m + 2}}}}}  - \frac{1}{2}\sum\limits_{n = 1}^\infty  {\frac{{{\zeta _n}\left( 2 \right)}}{{{n^{m + 1}}}}} ,\tag{3.7}\\
&\frac{1}{2}\sum\limits_{n = 1}^\infty  {\frac{{{}H_n^2}}{{{n^{m + 1}}}}}{\left( { - 1} \right)}^{n - 1}+ {\left( { - 1} \right)^{m - 1}}\sum\limits_{n = 1}^\infty  {\frac{{{}{H_n}{L_n}\left( 1 \right)}}{{{n^{m + 1}}}}}{\left( { - 1} \right)}^{n - 1}\\
&=\sum\limits_{j = 1}^{m - 1} {{{\left( { - 1} \right)}^{j - 1}}\bar \zeta \left( {m - j + 1} \right)\sum\limits_{n = 1}^\infty  {\frac{{{H_n}}}{{{n^{j + 1}}}}} }  + {\left( { - 1} \right)^{m - 1}}\ln 2\sum\limits_{n = 1}^\infty  {\frac{{{H_n}}}{{{n^{m + 1}}}}}  + \sum\limits_{n = 1}^\infty  {\frac{{{}H_n^{}}}{{{n^{m + 2}}}}}{\left( { - 1} \right)}^{n - 1} \\
&\quad + {\left( { - 1} \right)^{m - 1}}\ln 2\sum\limits_{n = 1}^\infty  {\frac{{{}{H_n}}}{{{n^{m + 1}}}}}{\left( { - 1} \right)}^{n - 1}- \frac{1}{2}\sum\limits_{n = 1}^\infty  {\frac{{{}{\zeta _n}\left( 2 \right)}}{{{n^{m + 1}}}}}{\left( { - 1} \right)}^{n - 1} - \zeta \left( 2 \right)\bar \zeta \left( {m + 1} \right),
 \tag{3.8}\\
 &\left( {\frac{1}{2} + {{\left( { - 1} \right)}^m}} \right)\sum\limits_{n = 1}^\infty  {\frac{{L_n^2\left( 1 \right)}}{{{n^{m + 1}}}}{{\left( { - 1} \right)}^{n - 1}}} \\
 & = \bar \zeta \left( 2 \right)\bar \zeta \left( {m + 1} \right) + \ln 2\sum\limits_{n = 1}^\infty  {\frac{{{H_n} + {L_n}\left( 1 \right)}}{{{n^{m + 1}}}}{{\left( { - 1} \right)}^{n - 1}}}  + {\left( { - 1} \right)^m}\ln 2\sum\limits_{n = 1}^\infty  {\frac{{{L_n}\left( 1 \right)}}{{{n^{m + 1}}}}\left( {1 + {{\left( { - 1} \right)}^{n - 1}}} \right)}\\
 &\quad  + \sum\limits_{n = 1}^\infty  {\frac{{{L_n}\left( 1 \right)}}{{{n^{m + 2}}}}}  - \ln 2\left( {\bar \zeta \left( {m + 2} \right) + \zeta \left( {m + 2} \right)} \right) - \frac{1}{2}\sum\limits_{n = 1}^\infty  {\frac{{{\zeta _n}\left( 2 \right)}}{{{n^{m + 1}}}}{{\left( { - 1} \right)}^{n - 1}}} \\
 &\quad - \sum\limits_{j = 1}^{m - 1} {{{\left( { - 1} \right)}^{j - 1}}} \bar \zeta \left( {m - j + 1} \right)\sum\limits_{n = 1}^\infty  {\frac{{{L_n}\left( 1 \right)}}{{{n^{j + 1}}}}} , \tag{3.9}\\
 &\frac{1}{2}\sum\limits_{n = 1}^\infty  {\frac{{L_n^2\left( 1 \right)}}{{{n^{m + 1}}}}}  + {\left( { - 1} \right)^{m - 1}}\sum\limits_{n = 1}^\infty  {\frac{{{H_n}{L_n}\left( 1 \right)}}{{{n^{m + 1}}}}}\\
 & = \bar \zeta \left( 2 \right)\zeta \left( {m + 1} \right) + \ln 2\sum\limits_{n = 1}^\infty  {\frac{{{H_n} + {L_n}\left( 1 \right)}}{{{n^{m + 1}}}}}  - \ln 2\left( {\bar \zeta \left( {m + 2} \right) + \zeta \left( {m + 2} \right)} \right) - \frac{1}{2}\sum\limits_{n = 1}^\infty  {\frac{{{\zeta _n}\left( 2 \right)}}{{{n^{m + 1}}}}} \\
 &\quad + \sum\limits_{n = 1}^\infty  {\frac{{{L_n}\left( 1 \right)}}{{{n^{m + 2}}}}} {\left( { - 1} \right)^{n - 1}} - \sum\limits_{j = 1}^{m - 1} {{{\left( { - 1} \right)}^{j - 1}}} \zeta \left( {m - j + 1} \right)\sum\limits_{n = 1}^\infty  {\frac{{{L_n}\left( 1 \right)}}{{{n^{j + 1}}}}} .\tag{3.10}
\end{align*}
\end{thm}
\pf Similarly as in the proofs of Theorem 3.1 and 3.2, we consider the following sums
\[\sum\limits_{k = 1}^\infty  {\sum\limits_{n = 1}^\infty  {\frac{{{H_n}}}{{{k^m}n\left( {n + k} \right)}}} } ,\sum\limits_{k = 1}^\infty  {\sum\limits_{n = 1}^\infty  {\frac{{{H_n}}}{{{k^m}n\left( {n + k} \right)}}{{\left( { - 1} \right)}^{k - 1}}} } ,\]\[\sum\limits_{k = 1}^\infty  {\sum\limits_{n = 1}^\infty  {\frac{{{L_n}\left( 1 \right)}}{{{k^m}n\left( {n + k} \right)}}} } ,\sum\limits_{k = 1}^\infty  {\sum\limits_{n = 1}^\infty  {\frac{{{L_n}\left( 1 \right)}}{{{k^m}n\left( {n + k} \right)}}} } {\left( { - 1} \right)^{k - 1}}.\]
Then using identities (2.12), (2.13) with the help of the following formula
\begin{align*}
\sum\limits_{n = 1}^\infty  {\frac{{{{\left( { - 1} \right)}^{n - 1}}}}{{{n^m}\left( {n + k} \right)}}}&= \sum\limits_{j = 1}^{m - 1} {\frac{{{{\left( { - 1} \right)}^{j - 1}}}}{{{k^j}}}\bar \zeta \left( {m - j + 1} \right)} \; + \frac{{{{\left( { - 1} \right)}^{m - 1}}}}{{{k^m}}}\ln 2 \\&\quad+ \frac{{{{\left( { - 1} \right)}^{m + k}}}}{{{k^m}}}\ln 2 - \frac{{{{\left( { - 1} \right)}^{m + k}}}}{{{k^m}}}{L_k}\left( 1 \right),
\end{align*}
we deduce Theorem 3.3 holds.\hfill$\square$
\begin{thm} For integer $m>0$, we have
\begin{align*}
&\left( {\frac{1}{3} + {{\left( { - 1} \right)}^m}} \right)\sum\limits_{n = 1}^\infty  {\frac{{H_n^3}}{{{n^{m + 1}}}}}  + \left( {1 + {{\left( { - 1} \right)}^{m - 1}}} \right)\sum\limits_{n = 1}^\infty  {\frac{{{H_n}{\zeta _n}\left( 2 \right)}}{{{n^{m + 1}}}}} \\
& =\sum\limits_{j = 0}^{m - 2} {{{\left( { - 1} \right)}^j}\zeta \left( {m - j} \right)} \sum\limits_{n = 1}^\infty  {\frac{{H_n^2 - {\zeta _n}\left( 2 \right)}}{{{n^{j + 2}}}}}  + \sum\limits_{n = 1}^\infty  {\frac{{H_n^2 + {\zeta _n}\left( 2 \right)}}{{{n^{m + 2}}}}}  - \frac{2}{3}\sum\limits_{n = 1}^\infty  {\frac{{{\zeta _n}\left( 3 \right)}}{{{n^{m + 1}}}}}  - 2\zeta \left( 3 \right)\zeta \left( {m + 1} \right),\tag{3.11}\\
&\sum\limits_{n = 1}^\infty  {\frac{{L_n^3\left( 1 \right) + {L_n}\left( 1 \right){\zeta _n}\left( 2 \right)}}{{{n^{2m + 1}}}}}\\
& = 2\bar \zeta \left( 2 \right)\left( {\sum\limits_{n = 1}^\infty  {\frac{{{L_n}\left( 1 \right)}}{{{n^{2m + 1}}}}} } \right) + 2\ln 2\sum\limits_{n = 1}^\infty  {\frac{{{H_n}{L_n}\left( 1 \right) + L_n^2\left( 1 \right)}}{{{n^{2m + 1}}}}}  - 2\ln 2\sum\limits_{n = 1}^\infty  {\frac{{{L_n}\left( 1 \right)}}{{{n^{2m + 2}}}}\left( {1 + {{\left( { - 1} \right)}^{n - 1}}} \right)}\\
&\quad+ 2\sum\limits_{n = 1}^\infty  {\frac{{L_n^2\left( 1 \right)}}{{{n^{2m + 2}}}}{{\left( { - 1} \right)}^{n - 1}}}  - 2\sum\limits_{i = 1}^m {{{\left( { - 1} \right)}^{i - 1}}} \left( {\sum\limits_{n = 1}^\infty  {\frac{{{L_n}\left( 1 \right)}}{{{n^{i + 1}}}}} } \right)\left( {\sum\limits_{n = 1}^\infty  {\frac{{{L_n}\left( 1 \right)}}{{{n^{2m + 1 - i}}}}} } \right)\\
&\quad + {\left( { - 1} \right)^{m - 1}}{\left( {\sum\limits_{n = 1}^\infty  {\frac{{{L_n}\left( 1 \right)}}{{{n^{m + 1}}}}} } \right)^2},\tag{3.12}\\
&\frac{1}{2}\sum\limits_{n = 1}^\infty  {\frac{{{H_n}L_n^2\left( 1 \right) + {H_n}{\zeta _n}\left( 2 \right)}}{{{n^{m + 1}}}}}  + \frac{{{{\left( { - 1} \right)}^{m - 1}}}}{2}\sum\limits_{n = 1}^\infty  {\frac{{H_n^2{L_n}\left( 1 \right) + {L_n}\left( 1 \right){\zeta _n}\left( 2 \right)}}{{{n^{m + 1}}}}}\\
&= \bar \zeta \left( 2 \right)\left( {\sum\limits_{n = 1}^\infty  {\frac{{{H_n}}}{{{n^{m + 1}}}}} } \right) + \ln 2\sum\limits_{n = 1}^\infty  {\frac{{H_n^2 + {H_n}{L_n}\left( 1 \right)}}{{{n^{m + 1}}}}}  - \ln 2\sum\limits_{n = 1}^\infty  {\frac{{{H_n}}}{{{n^{m + 2}}}}\left( {1 + {{\left( { - 1} \right)}^{n - 1}}} \right)}  \\
&\quad - \sum\limits_{i = 1}^{m - 1} {{{\left( { - 1} \right)}^{i - 1}}} \left( {\sum\limits_{n = 1}^\infty  {\frac{{{L_n}\left( 1 \right)}}{{{n^{i + 1}}}}} } \right)\left( {\sum\limits_{n = 1}^\infty  {\frac{{{H_n}}}{{{n^{m + 1 - i}}}}} } \right) - {\left( { - 1} \right)^{m - 1}}\zeta \left( 2 \right)\left( {\sum\limits_{n = 1}^\infty  {\frac{{{L_n}\left( 1 \right)}}{{{n^{m + 1}}}}} } \right) \\
&\quad+ {\left( { - 1} \right)^{m - 1}}\sum\limits_{n = 1}^\infty  {\frac{{{H_n}{L_n}\left( 1 \right)}}{{{n^{m + 2}}}}} + \sum\limits_{n = 1}^\infty  {\frac{{{H_n}{L_n}\left( 1 \right)}}{{{n^{m + 2}}}}{{\left( { - 1} \right)}^{n - 1}}}.\tag{3.13}
\end{align*}
\end{thm}
\pf Similarly as in the proof of Theorem 3.1-3.3, we consider the following sums
\[\sum\limits_{k = 1}^\infty  {\sum\limits_{n = 1}^\infty  {\frac{{H_n^2 - {\zeta _n}\left( 2 \right)}}{{{k^m}n\left( {n + k} \right)}}} } ,\sum\limits_{k = 1}^\infty \sum\limits_{n = 1}^\infty  {\frac{{{L_k}\left( 1 \right){L_n}\left( 1 \right)}}{{{k^{2m}}n\left( {n + k} \right)}}} ,\sum\limits_{k = 1}^\infty \sum\limits_{n = 1}^\infty  {\frac{{{H_k}{L_n}\left( 1 \right)}}{{{k^m}n\left( {n + k} \right)}}} .\]
Then using (2.12), (2.13) and (2.28), by a simple calculation, we obtain the desired results.\hfill$\square$ \\
Letting $p_1=p_2=p,m=2k$ in Theorem 3.2, we can give the following Corollary
\begin{cor}
For integers $k>0,p>1$, we have
\begin{align*}
&\frac{{\left( {p - 1} \right)!}}{p}\sum\limits_{n = p - 1}^\infty  {\frac{{S\left( {n + 1,p} \right){Y_p}\left( n \right)}}{{{n^{2k + 1}}n!}}} \\
& = \left( {p - 1} \right)!\sum\limits_{n = p - 1}^\infty  {\frac{{S\left( {n + 1,p} \right){Y_{p - 1}}\left( n \right)}}{{{n^{2k + 2}}n!}}}  - {\left[ {\left( {p - 1} \right)!} \right]^2}\zeta \left( p \right)\left( {\sum\limits_{n = p - 1}^\infty  {\frac{{S\left( {n + 1,p} \right)}}{{{n^{2k + 1}}n!}}} } \right)\\
&\quad + {\left[ {\left( {p - 1} \right)!} \right]^2}\sum\limits_{i = 1}^k {{{\left( { - 1} \right)}^{i - 1}}\left( {\sum\limits_{n = p - 1}^\infty  {\frac{{S\left( {n + 1,p} \right)}}{{{n^{2k + 1 - i}}n!}}} } \right)\left( {\sum\limits_{n = p - 1}^\infty  {\frac{{S\left( {n + 1,p} \right)}}{{{n^{i + 1}}n!}}} } \right)} \\
&\quad - \frac{{{{\left[ {\left( {p - 1} \right)!} \right]}^2}}}{2}{\left( { - 1} \right)^{k - 1}}{\left( {\sum\limits_{n = p - 1}^\infty  {\frac{{S\left( {n + 1,p} \right)}}{{{n^{k + 1}}n!}}} } \right)^2}.\tag{3.14}
\end{align*}
\end{cor}
Taking $p=2$ in (3.14), we obtain
\begin{align*}
\sum\limits_{n = 1}^\infty  {\frac{{H_n^3 + {H_n}{\zeta _n}\left( 2 \right)}}{{{n^{2k + 1}}}}}  &= 2\sum\limits_{n = 1}^\infty  {\frac{{H_n^2}}{{{n^{2k + 2}}}}}  + 2\sum\limits_{i = 1}^k {{{\left( { - 1} \right)}^{i - 1}}\left( {\sum\limits_{n = 1}^\infty  {\frac{{{H_n}}}{{{n^{2k + 1 - i}}}}} } \right)} \left( {\sum\limits_{n = 1}^\infty  {\frac{{{H_n}}}{{{n^{i + 1}}}}} } \right)\\
 &\quad- 2\zeta \left( 2 \right)\left( {\sum\limits_{n = 1}^\infty  {\frac{{{H_n}}}{{{n^{2k + 1}}}}} } \right) - {\left( { - 1} \right)^{k - 1}}{\left( {\sum\limits_{n = 1}^\infty  {\frac{{{H_n}}}{{{n^{k + 1}}}}} } \right)^2}.\tag{3.15}
\end{align*}
Letting $m=2k-1$ in (3.11), we get
\begin{align*}
\sum\limits_{n = 1}^\infty  {\frac{{H_n^3}}{{{n^{2k + 1}}}}}
 &=\frac{3}{4}\sum\limits_{j = 0}^{2k - 2} {{{\left( { - 1} \right)}^j}\zeta \left( {2k - j} \right)} \sum\limits_{n = 1}^\infty  {\frac{{H_n^2 - {\zeta _n}\left( 2 \right)}}{{{n^{j + 2}}}}}  + \frac{3}{4}\sum\limits_{n = 1}^\infty  {\frac{{H_n^2 + {\zeta _n}\left( 2 \right)}}{{{n^{2k + 2}}}}}
  \nonumber \\
           &\quad - \frac{1}{2}\sum\limits_{n = 1}^\infty  {\frac{{{\zeta _n}\left( 3 \right)}}{{{n^{2k + 1}}}}}  - \frac{3}{2}\zeta \left( 3 \right)\zeta \left( {2k + 1} \right).\tag{3.16}
\end{align*}
Substituting (3.16) into (3.15), we arrive at the conclusion that
\begin{align*}
\sum\limits_{n = 1}^\infty  {\frac{{{H_n}{\zeta _n}\left( 2 \right)}}{{{n^{2k + 1}}}}}
 &=2\sum\limits_{n = 1}^\infty  {\frac{{H_n^2}}{{{n^{2k + 2}}}}}  + 2\sum\limits_{i = 1}^k {{{\left( { - 1} \right)}^{i - 1}}\left( {\sum\limits_{n = 1}^\infty  {\frac{{{H_n}}}{{{n^{2k + 1 - i}}}}} } \right)} \left( {\sum\limits_{n = 1}^\infty  {\frac{{{H_n}}}{{{n^{i + 1}}}}} } \right)
  \nonumber \\
           &\quad + \frac{1}{2}\sum\limits_{n = 1}^\infty  {\frac{{{\zeta _n}\left( 3 \right)}}{{{n^{2k + 1}}}}}  + \frac{3}{2}\zeta \left( 3 \right)\zeta \left( {2k + 1} \right) - 2\zeta \left( 2 \right)\left( {\sum\limits_{n = 1}^\infty  {\frac{{{H_n}}}{{{n^{2k + 1}}}}} } \right) - {\left( { - 1} \right)^{k - 1}}{\left( {\sum\limits_{n = 1}^\infty  {\frac{{{H_n}}}{{{n^{k + 1}}}}} } \right)^2}
  \nonumber \\
           &\quad - \frac{3}{4}\sum\limits_{j = 0}^{2k - 2} {{{\left( { - 1} \right)}^j}\zeta \left( {2k - j} \right)} \sum\limits_{n = 1}^\infty  {\frac{{H_n^2 - {\zeta _n}\left( 2 \right)}}{{{n^{j + 2}}}}}  - \frac{3}{4}\sum\limits_{n = 1}^\infty  {\frac{{H_n^2 + {\zeta _n}\left( 2 \right)}}{{{n^{2k + 2}}}}} .\tag{3.17}
\end{align*}
Similarly, taking $(p_1,p_2)=(2,3),(1,4)$ in Theorem 3.2, we deduce that
\begin{align*}
&\left( {\frac{1}{3} - \frac{{{{\left( { - 1} \right)}^{m - 1}}}}{2}} \right)\sum\limits_{n = 1}^\infty  {\frac{{H_n^4}}{{{n^{m + 1}}}}}  + \sum\limits_{n = 1}^\infty  {\frac{{H_n^2{\zeta _n}\left( 2 \right)}}{{{n^{m + 1}}}}}  + \frac{2}{3}\sum\limits_{n = 1}^\infty  {\frac{{{H_n}{\zeta _n}\left( 3 \right)}}{{{n^{m + 1}}}}}  + \frac{{{{\left( { - 1} \right)}^{m - 1}}}}{2}\sum\limits_{n = 1}^\infty  {\frac{{\zeta _n^2\left( 2 \right)}}{{{n^{m + 1}}}}} \\
& = \left( {1 - {{\left( { - 1} \right)}^{m - 1}}} \right)\sum\limits_{n = 1}^\infty  {\frac{{H_n^3}}{{{n^{m + 2}}}}}  + \left( {1 + {{\left( { - 1} \right)}^{m - 1}}} \right)\sum\limits_{n = 1}^\infty  {\frac{{{H_n}{\zeta _n}\left( 2 \right)}}{{{n^{m + 2}}}}}  + {\left( { - 1} \right)^{m - 1}}\zeta \left( 2 \right)\sum\limits_{n = 1}^\infty  {\frac{{H_n^2 - {\zeta _n}\left( 2 \right)}}{{{n^{m + 1}}}}}\\
&\quad - 2\zeta \left( 3 \right)\sum\limits_{n = 1}^\infty  {\frac{{{H_n}}}{{{n^{m + 1}}}}}  + \sum\limits_{i = 1}^{m - 1} {{{\left( { - 1} \right)}^{i - 1}}} \left( {\sum\limits_{n = 1}^\infty  {\frac{{{H_n}}}{{{n^{m + 1 - i}}}}} } \right)\left( {\sum\limits_{n = 1}^\infty  {\frac{{H_n^2 - {\zeta _n}\left( 2 \right)}}{{{n^{i + 1}}}}} } \right).\tag{3.18}
\end{align*}
\begin{align*}
&\left( {\frac{1}{4} + {{\left( { - 1} \right)}^m}} \right)\sum\limits_{n = 1}^\infty  {\frac{{H_n^4}}{{{n^{m + 1}}}}}  + 3\left( {\frac{1}{2} + {{\left( { - 1} \right)}^{m - 1}}} \right)\sum\limits_{n = 1}^\infty  {\frac{{H_n^2{\zeta _n}\left( 2 \right)}}{{{n^{m + 1}}}}}  + 2\left( {1 + {{\left( { - 1} \right)}^m}} \right)\sum\limits_{n = 1}^\infty  {\frac{{{H_n}{\zeta _n}\left( 3 \right)}}{{{n^{m + 1}}}}}  + \frac{3}{4}\sum\limits_{n = 1}^\infty  {\frac{{\zeta _n^2\left( 2 \right)}}{{{n^{m + 1}}}}}\\
&=\sum\limits_{n = 1}^\infty  {\frac{{H_n^3 + 3{H_n}{\zeta _n}\left( 2 \right) + 2{\zeta _n}\left( 3 \right)}}{{{n^{m + 2}}}}}  + \sum\limits_{i = 1}^{m - 1} {{{\left( { - 1} \right)}^{i - 1}}} \zeta \left( {m + 1 - i} \right)\left( {\sum\limits_{n = 1}^\infty  {\frac{{H_n^3 - 3{H_n}{\zeta _n}\left( 2 \right) + 2{\zeta _n}\left( 3 \right)}}{{{n^{i + 1}}}}} } \right)\\
&\quad - \frac{3}{2}\sum\limits_{n = 1}^\infty  {\frac{{{\zeta _n}\left( 4 \right)}}{{{n^{m + 1}}}}}  - 6\zeta \left( 4 \right)\zeta \left( {m + 1} \right).\tag{3.19}
\end{align*}
Proceeding in a similar fashion to evaluation of the Theorem 3.1-3.4, it is possible to evaluate other Euler sums involving harmonic numbers and alternating harmonic numbers. For instance, multiplying (2.21) by $\frac{{{{\left( { - 1} \right)}^{k - 1}}}}{{{k^m}}}, \frac{{{L_k}\left( 1 \right)}}{{{k^m}}}$ and summing with respect to $k$, we obtain
\begin{align*}
&\frac{1}{p}\sum\limits_{n = 1}^\infty  {\frac{{{Y_p}\left( n \right)}}{{{n^{m + 1}}}}{{\left( { - 1} \right)}^{n - 1}}}  + {\left( { - 1} \right)^{m - 1}}\left( {p - 1} \right)!\sum\limits_{n = p - 1}^\infty  {\frac{{S\left( {n + 1,p} \right){L_n}\left( 1 \right)}}{{n!{n^{m + 1}}}}{{\left( { - 1} \right)}^{n - 1}}} \\
&= \sum\limits_{n = 1}^\infty  {\frac{{{Y_{p - 1}}\left( n \right)}}{{{n^{m + 2}}}}{{\left( { - 1} \right)}^{n - 1}}}  + \left( {p - 1} \right)!\sum\limits_{i = 1}^{m - 1} {{{\left( { - 1} \right)}^{i - 1}}\bar \zeta \left( {m + 1 - i} \right)\left( {\sum\limits_{n = p - 1}^\infty  {\frac{{S\left( {n + 1,p} \right)}}{{n!{n^{i + 1}}}}} } \right)} \\
&\quad + {\left( { - 1} \right)^{m - 1}}\left( {p - 1} \right)!\ln 2\sum\limits_{n = p - 1}^\infty  {\frac{{S\left( {n + 1,p} \right)}}{{n!{n^{m + 1}}}}\left( {1 + {{\left( { - 1} \right)}^{n - 1}}} \right)}  - \left( {p - 1} \right)!\zeta \left( p \right)\bar \zeta \left( {m + 1} \right),\tag{3.20}
\end{align*}
and
\begin{align*}
&\frac{1}{p}\sum\limits_{n = 1}^\infty  {\frac{{{Y_p}\left( n \right){L_n}\left( 1 \right)}}{{{n^{m + 1}}}}}  + \frac{{{{\left( { - 1} \right)}^{m - 1}}}}{2}\left( {p - 1} \right)!\sum\limits_{n = p - 1}^\infty  {\frac{{S\left( {n + 1,p} \right)\left( {L_n^2\left( 1 \right) + {\zeta _n}\left( 2 \right)} \right)}}{{n!{n^{m + 1}}}}} \\
&=\sum\limits_{n = 1}^\infty  {\frac{{{Y_{p - 1}}\left( n \right){L_n}\left( 1 \right)}}{{{n^{m + 2}}}}}  + \left( {p - 1} \right)!\sum\limits_{i = 1}^{m - 1} {{{\left( { - 1} \right)}^{i - 1}}\left( {\sum\limits_{n = 1}^\infty  {\frac{{{L_n}\left( 1 \right)}}{{{n^{m + 1 - i}}}}} } \right)\left( {\sum\limits_{n = p - 1}^\infty  {\frac{{S\left( {n + 1,p} \right)}}{{n!{n^{i + 1}}}}} } \right)}\\
&\quad+ {\left( { - 1} \right)^{m - 1}}\left( {p - 1} \right)!\bar \zeta \left( 2 \right)\sum\limits_{n = p - 1}^\infty  {\frac{{S\left( {n + 1,p} \right)}}{{n!{n^{m + 1}}}}} + {\left( { - 1} \right)^{m - 1}}\left( {p - 1} \right)!\sum\limits_{n = p - 1}^\infty  {\frac{{S\left( {n + 1,p} \right){L_n}\left( 1 \right)}}{{n!{n^{m + 2}}}}} \\
&\quad - {\left( { - 1} \right)^{m - 1}}\left( {p - 1} \right)!\ln 2\sum\limits_{n = p - 1}^\infty  {\frac{{S\left( {n + 1,p} \right)}}{{n!{n^{m + 2}}}}\left( {1 + {{\left( { - 1} \right)}^{n - 1}}} \right)}  \\
&\quad + {\left( { - 1} \right)^{m - 1}}\left( {p - 1} \right)!\ln 2\sum\limits_{n = p - 1}^\infty  {\frac{{S\left( {n + 1,p} \right)\left( {{H_n} + {L_n}\left( 1 \right)} \right)}}{{n!{n^{m + 1}}}}}\\
&\quad- \left( {p - 1} \right)!\zeta \left( p \right)\sum\limits_{n = 1}^\infty  {\frac{{{L_n}\left( 1 \right)}}{{{n^{m + 1}}}}}.\tag{3.21}
\end{align*}
\section{Some Examples}
Now, we give some examples.
\begin{align*}
&\sum\limits_{n = 1}^\infty  {\frac{{L_n^2\left( 1 \right)}}{{{n^2}}}} {\left( { - 1} \right)^{n - 1}} =  - \frac{{41}}{{16}}\zeta \left( 4 \right) + 2\zeta \left( 2 \right){\ln ^2}2 + \frac{1}{6}{\ln ^4}2 + \frac{7}{4}\zeta \left( 3 \right)\ln 2 + 4{\rm{L}}{{\rm{i}}_4}\left( {\frac{1}{2}} \right),\\
&\sum\limits_{n = 1}^\infty  {\frac{{L_n^2\left( 1 \right)}}{{{n^3}}}} {\left( { - 1} \right)^{n - 1}} =  - 4{\rm{L}}{{\rm{i}}_4}\left( {\frac{1}{2}} \right)\ln 2 + \frac{{19}}{8}\zeta \left( 4 \right)\ln 2 + \zeta \left( 2 \right){\ln ^3}2 - \frac{1}{6}{\ln ^5}2 + \frac{3}{8}\zeta \left( 2 \right)\zeta \left( 3 \right) - \frac{{19}}{{32}}\zeta \left( 5 \right),\\
&\sum\limits_{n = 1}^\infty  {\frac{{L_n^2\left( 1 \right)}}{{{n^4}}}{{\left( { - 1} \right)}^{n - 1}}}  = \frac{{15}}{4}{\ln ^2}2\zeta \left( 4 \right) + \frac{9}{4}\zeta \left( 2 \right)\zeta \left( 3 \right)\ln 2 - \frac{{93}}{{16}}\zeta \left( 5 \right)\ln 2 + \frac{{35}}{{64}}\zeta \left( 6 \right) - \frac{{15}}{{16}}{\zeta ^2}\left( 3 \right) \\
&\quad \quad\quad\quad\quad\quad\quad \quad\quad+ \sum\limits_{n = 1}^\infty  {\frac{{{\zeta _n}\left( 2 \right)}}{{{n^4}}}{{\left( { - 1} \right)}^{n - 1}}} , \\
&\sum\limits_{n = 1}^\infty  {\frac{{{L_n}\left( 1 \right){L_n}\left( 2 \right)}}{n}} {\left( { - 1} \right)}^{n - 1} = \frac{{61}}{{16}}\zeta \left( 4 \right) - \frac{7}{8}\zeta \left( 3 \right)\ln 2 - \frac{1}{4}\zeta \left( 2 \right){\ln ^2}2 - \frac{1}{6}{\ln ^4}2 - 4{\rm{L}}{{\rm{i}}_4}\left( {\frac{1}{2}} \right),\\
&\sum\limits_{n = 1}^\infty  {\frac{{{L_n}\left( 1 \right){L_n}\left( 3 \right)}}{n}} {\left( { - 1} \right)^{n - 1}} = 2\ln 2{\rm{L}}{{\rm{i}}_4}\left( {\frac{1}{2}} \right) + \frac{1}{{12}}{\ln ^5}2 + \frac{3}{8}\zeta \left( 3 \right){\ln ^2}2 - \frac{{19}}{{32}}\zeta \left( 5 \right) - \frac{1}{2}\zeta \left( 2 \right){\ln ^3}2\\
&\quad\quad \quad\quad\quad\quad\quad \quad\quad\quad\quad\quad  + \frac{{11}}{{16}}\zeta \left( 4 \right)\ln 2 + \frac{1}{4}\zeta \left( 2 \right)\zeta \left( 3 \right),\\
&\sum\limits_{n = 1}^\infty  {\frac{{{L_n}\left( 1 \right){L_n}\left( 4 \right)}}{n}} {\left( { - 1} \right)^{n - 1}} =  - \frac{{35}}{{128}}\zeta \left( 6 \right) + \frac{3}{4}{\zeta ^2}\left( 3 \right) - \frac{9}{8}\zeta \left( 2 \right)\zeta \left( 3 \right)\ln 2 + \frac{{155}}{{32}}\zeta \left( 5 \right)\ln 2\\
&\quad\quad \quad\quad\quad\quad\quad \quad\quad\quad\quad\quad \quad - \frac{{23}}{{16}}\zeta \left( 4 \right){\ln ^2}2 - \sum\limits_{n = 1}^\infty  {\frac{{{\zeta _n}\left( 2 \right)}}{{{n^4}}}{{\left( { - 1} \right)}^{n - 1}}} ,\\
&\sum\limits_{n{\rm{ = }}1}^\infty  {\frac{{H_n^3}}{{{n^5}}}}  = \frac{{469}}{{32}}\zeta \left( 8 \right) - 16\zeta \left( 3 \right)\zeta \left( 5 \right) + \frac{3}{2}\zeta \left( 2 \right){\zeta ^2}\left( 3 \right) + \frac{{11}}{4}\sum\limits_{n = 1}^\infty  {\frac{{{\zeta _n}\left( 2 \right)}}{{{n^6}}}},\\
& \sum\limits_{n = 1}^\infty  {\frac{{{H_n}{\zeta _n}\left( 2 \right)}}{{{n^5}}}}  =  - \frac{{343}}{{48}}\zeta \left( 8 \right) + 12\zeta \left( 3 \right)\zeta \left( 5 \right) - \frac{5}{2}\zeta \left( 2 \right){\zeta ^2}\left( 3 \right) - \frac{3}{4}\sum\limits_{n = 1}^\infty  {\frac{{{\zeta _n}\left( 2 \right)}}{{{n^6}}}},\\
& \sum\limits_{n = 1}^\infty  {\frac{{H_n^2{\zeta _n}\left( 3 \right)}}{{{n^3}}}}  = \frac{9}{2}\zeta \left( 3 \right)\zeta \left( 5 \right) + \frac{3}{2}\zeta \left( 2 \right){\zeta ^2}\left( 3 \right) - \frac{{443}}{{288}}\zeta \left( 8 \right) - \frac{23}{4}\sum\limits_{n = 1}^\infty  {\frac{{{\zeta _n}\left( 2 \right)}}{{{n^6}}}}.
\end{align*}
{\bf Acknowledgments.} The authors would like to thank the anonymous
referee for his/her helpful comments, which improve the presentation
of the paper.
\end{document}